\documentclass[11pt]{article}
\usepackage{amssymb}
\usepackage{amsmath} 
\usepackage{amsfonts}
\usepackage{amsthm}
\usepackage{graphicx}
\usepackage{subcaption}
\usepackage{lmodern}
\usepackage{color}
\usepackage{comment}

\newtheorem{theorem}{Theorem}[section]
\newtheorem{lemma}[theorem]{Lemma}
\newtheorem{corollary}[theorem]{Corollary}

\newtheorem{conjecture}[theorem]{Conjecture}

\newtheorem{remark}[theorem]{Remark}
\usepackage{fullpage}
\usepackage{authblk}

\newcommand{\ud}{{\rm d}}
\newcommand{\error}{e_{\varepsilon}}

\begin{document}

\title{Space-time finite element methods \\ for nonlinear wave equations \\ via elliptic regularisation}
\author[1]{Lehel Banjai}
\author[1,2,3]{Emmanuil H.~Georgoulis}
\author[1]{Brian Hennessy\thanks{bgh2000@hw.ac.uk}}
\affil[1]{The Maxwell Institute for Mathematical Sciences \& Department of Mathematics, School of Mathematical and Computer Sciences, Heriot-Watt University, Edinburgh, UK}
\affil[2]{ Department of Mathematics, School of Applied Mathematics and Physical Sciences, National Technical University of Athens, Zografou Campus, Athens, Greece}
\affil[3]{ Institute for Applied and Computational Mathematics-FORTH, Heraklion-Crete, Greece}
\maketitle
\begin{abstract}
\noindent We present and analyse a new conforming space-time Galerkin discretisation of a semi-linear wave equation, based on a variational formulation derived from De Giorgi's elliptic regularisation viewpoint of the wave equation in second-order formulation. The method is shown to be well-posed through a minimisation approach, and also unconditionally stable for all choices of conforming discretisation spaces. Further, \emph{a priori} error bounds are proven for sufficiently smooth solutions. Special attention is given to the conditioning of the method and its stable implementation. Numerical experiments are provided to validate the theoretical findings.
\end{abstract}
\section{Introduction}
Nonlinear wave equations are abundant in the mathematical modelling of complex wave phenomena, in mathematical relativity, and in quantum mechanics. Many such equations are characterised by dispersion effects, posing significant challenges for the numerical approximation of their solutions. Further challenges include the preservation of important properties of exact solutions at the discrete level or, even the well-posedness of numerical solutions themselves in some instances. Finite element methods are of interest to allow for local high resolution in space-time regions where needed, or approximation on complex bounded spatial domains. 

Space-time Galerkin/finite element methods have enjoyed recent interest in the literature, as they enable highly localised resolution in space-time. This is of particular interest in nonlinear wave equations whereby the evolution may result to, typically \emph{a priori} unknown localised frequencies and amplitudes in the solution. Also, local or global high wave velocities may hinder the performance of standard, explicit timestepping methods, rendering space-time methods a potentially attractive alternative.  

For concreteness, we consider the defocusing semilinear wave equation on $D:=\Omega\times I$, $I:=(0,T]$, $T>0$ (finite) final time, where $\Omega\subset \mathbb{R}^d$, $d = 1,2,3$, is a bounded, polygonal/polyhedral domain, reading: find $u: D\to\mathbb{R}$, such that 
\begin{equation}
\begin{aligned}
    u_{tt}-\Delta u + \frac{p}{2}|u|^{p-2}u&=f, & &\text{in } \Omega \times I \\
    u &= 0,& &\text{on }\partial\Omega \times I\\
    u(x,0)&=u^0(x), & &\text{in }\Omega\\
    u_t(x,0)&=u^1(x),&  &\text{in }\Omega,
\end{aligned}
\label{eq: Model equation}
\end{equation}
where $f \in L^2(0,T;L^2(\Omega))$ is a given source term, $p=0,2$ or $p>2$ which corresponds to the linear wave equation, the linear Klein-Gordon equation and the semilinear wave equation respectively, $d = 1,2,3$, and  $u^0\in H^1(\Omega),\;u^1\in L^2(\Omega)$, using standard Lebesgue and Sobolev space notation. This equation is used as a prototypical model for stable nonlinear dispersion and appears in many application areas, e.g., in models for nonlinear stability of light rings around compact stellar objects \cite{MR4921049,hy3r-ww3w}.
When $\Omega=\mathbb{R}^d$, 
local well-posedness is characterised by the nonlinearity $p$ and the dimension $d$. For $d > 1$, the equation is $H^1$-subcritical if $p < p_c = \frac{2d}{d-2}$, $H^1$-critical if $p= p_c$ and $H^1$-supercritical if $p >p_c$. In the supercritical regime, one expects that the high frequencies are unstable and develop nonlinear behaviour quickly so well-posedness is not available \cite[Principle~3.1]{MR2233925}. In dimension $d=1,2$ local well-posedness may be established using the Sobolev embeddings for any $p< \infty$. This also holds in $d=3$ up until $p=4$. For higher powers $4<p\leq 6$, one must exploit the dispersive character of the wave equation through the Strichartz estimates \cite[Theorem~5.2]{MR2233925}. On bounded domains with smooth boundary, well-posedness is still available but the Strichartz estimates require stronger admissibility criteria compared to the unbounded domain $\mathbb{R}^3$ due to microlocal phenomena such as caustics generated
in arbitrarily small time near the boundary \cite{MR2993000,MR2393429,MR2966654}.
For the defocusing problem, this regularity for the initial data is of interest as it is required to promote the locally well-posed problem to a globally well-posed one \cite[Section~3.4]{MR2233925}.

The above brief (and incomplete) account of available well-posedness results for the model problem \eqref{eq: Model equation} highlights the respective challenges that arise in the proof of rigorous results in numerical analysis for these problems. Variational approaches in the analysis of nonlinear wave equations naturally give rise to numerical methods of the Galerkin or Petrov-Galerkin type.
Space-time finite element methods for linear wave problems have appeared since the late 1980s; we refer to 
\cite{MR1241481,MR928689,MR1241479} for some early pioneering works. Their development has continued to date with discontinuous Galerkin approaches \cite{MR3592080,MR4286256,MR3518361,MR2272600,MR4550324,MR3556398,MR3748308,MR4083033} and methods based on inf-sup stability results and/or on Petrov-Galerkin approaches  \cite{MR4054692,MR4102914}.

Also, unconditionally stable methods based on second-order formulation for the linear wave equation have been proposed, e.g., by adding non-consistent penalty terms into the variational formulation. In \cite{MR4054692} a stabilization based on piecewise linear continuous functions is studied. This approach has been generalised for higher-order continuous polynomials and high order maximal regularity splines, see \cite{ferrari2025stability,FRASCHINI2024205,Z2021}. Related methods include using trial-to-test transformation operators \cite{ferrari2025space,MR4102914}.  Other variational approaches include ultraweak formulations, where integration by parts is used to move all derivatives onto the test function \cite{MR4444531}.  First-order formulations of the linear wave equation, while doubling the number of independent variables, can be unconditionally stable without the need for stabilization terms \cite{MR1255008,ferrari2024unconditionally,MR1325867,gomez2025variational,MR3556398,MR3748308}. 

To the best of our knowledge, none of these space-time approaches has been extended to second-order nonlinear wave problems. Apart from the attractive attribute of accepting locally variable numerical resolution, space-time methods enjoy a number of desirable properties, including high-order approximation in both space and time, the ability to use unstructured space-time meshes, multilevel preconditioning, and space–time parallelization. These attributes of space-time approaches become more pertinent when moving to nonlinear problems, due to the wealth of potential localised solution features that may arise during computation.

In this work, we propose and analyse a method for the second order in time formulation of the semi-linear wave equation \eqref{eq: Model equation} based on the elliptic regularisation in time. The methodology is motivated by De Giorgi's elliptic regularisation conjecture for hyperbolic problems, proposed in \cite{de1996congetture}. The basis idea in De Giorgi's approach is the definition of a nonstandard, regularisation-parameter-dependent family of \emph{convex} energies whose minimisers approximate the solution to \eqref{eq: Model equation} as the regularisation parameter $\varepsilon\to 0$. This approach naturally leads to a variational formulation for the wave equation which is elliptic, allowing for a well-posed space-time Galerkin discretisation to be defined. Owing to the ellipticity, an unconditionally stable second-order formulation for the nonlinear wave equation is obtained. The Galerkin setting allows eventually to show the existence and uniqueness of the discrete solution, as well as for quasi-optimality and \emph{a priori} error analysis for sufficiently smooth solutions in suitable exponentially weighted norms.  The results below also extend with minor modifications to other equations, e.g., the semilinear Klein-Gordon equation, (corresponding to \eqref{eq: Model equation} with a reaction term of the form $u+\frac{p}{2}|u|^{p-2}u$), which corresponds to self-interactions in Quantum Field Theory when $p=4$ \cite[Chapter~4]{MR1402248}. The effect of the regularisation parameter $\varepsilon$ is central to the study and its effects in conditioning are detailed. In particular, we show that care should be taken to alleviate numerical stability issues related to exponentially small numbers. A series of numerical experiments assessing the practicality of the methodology is also presented. 

The starting point of the present work is the exploration of the use of elliptic reguralisation in the context of rigorous numerical methods for nonlinear wave equation. We believe that, understanding the extent to which this analytic framework can motivate the construction of respective Galerkin methods, is ultimately of interest in the context of highly nonlinear complex wave problems that lack satisfactory numerical treatment. This work is, therefore, a first step in this direction. 

The remainder of this work is structured as follows. 
In Section \ref{sec:model_problem} we give a brief account of the elliptic reguralisation framework of De Giorgi and its relevance on the existence and uniqueness theory of nonlinear wave equations. In Section \ref{sec: discrete min} we detail the space-time finite element discretisation of the problem and prove well-posedness of the numerical method. In Section \ref{sec: Error rates}, we derive a quasi-optimality estimate and, consequently, convergence rates based on the regularity of the elliptic-regularised problem. In Section \ref{sec: numerical experiments}, we present some challenges in the implementation of the scheme, as well as numerical results as validation of the theoretical analysis. Finally in Section \ref{sec:conclusions} we draw some conclusions. 

\section{Elliptic regularisation for the wave equation}\label{sec:model_problem}
Elliptic regularisation for hyperbolic problems was conjectured by De Giorgi in \cite{de1996congetture}; the statement is taken from the English translation in \cite{MR3185411}.
\begin{conjecture}[De Giorgi \cite{de1996congetture}] Let $u^0,u^1 \in C^{\infty}_0(\mathbb{R}^d)$ and $p>2$ be an even integer. For every $\varepsilon> 0$, let $u^{\varepsilon}$ be the unique minimiser of the functional for $T=+\infty$ 
\begin{equation}
    I_{\varepsilon}(u):=\int_0^T\int_{\mathbb{R}^d}e^{-\frac{t}{\varepsilon}}\Big(|u_{tt}|^2+\frac{1}{\varepsilon^2}|\nabla u|^2+\frac{1}{\varepsilon^2}|u|^p\Big)\,\ud x\,\ud t    \label{eq: continuous functional no f}
\end{equation}
in the class of all functions satisfying $u^{\varepsilon}(x,0)=u^0,u^{\varepsilon}_t(x,0)=u^1$. 
Then there exists $\lim\limits_{\varepsilon\rightarrow 0}u^{\varepsilon}(x,t)=u(x,t)$, and satisfies the equation
\begin{equation}
    u_{tt}-\Delta u + \frac{p}{2}|u|^{p-2}u = 0.
\label{eq: semilin wave equation}
\end{equation}
\label{eq: De Giorgi Conjecture}
\end{conjecture}
The advantage of this approach is that it relates the solution of a non-convex problem, the semilinear wave equation \eqref{eq: semilin wave equation} to the minimisation problem of the convex functional \eqref{eq: continuous functional no f}. In addition, in supercritical cases where no uniqueness is guaranteed for \eqref{eq: semilin wave equation}, the functional \emph{still} admits a unique minimiser.

To see the connection of the minimisation problem to the solution of the semilinear wave equation, we formally compute the Euler-Lagrange equation of \eqref{eq: continuous functional no f}, which is given by
\begin{equation}
    \varepsilon^2u_{tttt}^{\varepsilon}-2\varepsilon u_{ttt}^{\varepsilon}+u_{tt}^{\varepsilon}-\Delta u^{\varepsilon}+\frac{p}{2}|u^{\varepsilon}|^{p-2}u^{\varepsilon}=0.
    \label{eq: fourth order pde}
\end{equation}
Setting formally $\varepsilon=0$, we recover the nonlinear wave equation \eqref{eq: semilin wave equation}. Note that \eqref{eq: fourth order pde} is an equation of elliptic type.

The first positive resolution of De Giorgi's conjecture was given in 2011 by Stefanelli in \cite{MR2819200}, for a slightly altered statement, assuming in particular a finite time horizon $T< \infty$, and no restriction on $p$ being even or an integer. A key element of Stefanelli's proof is the use of a finite difference discretisation of the time derivatives, along with an appropriate approximation of the exponential weight in \eqref{eq: continuous functional no f}. An $\varepsilon$-independent energy bound is established on the discrete level, which when followed by a $\Gamma$-convergence argument to take the finite difference discretisation to zero, an $\varepsilon$-independent energy bound on the minimisers $u^{\varepsilon}$ of \eqref{eq: continuous functional no f} is then established, which allows the passage of $\varepsilon\rightarrow 0$.

While existence and uniqueness to the set of discrete solutions can be guaranteed, the underlying method converges at first order in time. 
For the infinite time horizon case ($T=+\infty$), the conjecture was proved in 2012 by Serra and Tilli \cite{MR2912711}. Their argument is purely variational, making no use of the Euler-Lagrange equation \eqref{eq: fourth order pde}.

This type of variational regularization falls under the Weighted Inertia-Energy-Dissipation (WIDE) principle: this is a
global variational approach to nonlinear parabolic and hyperbolic equations, whereby minimisation of the $\varepsilon$-dependent WIDE functional on trajectories delivers an elliptic-in-time regularization. This approach has been applied to the study of highly nonlinear evolution problems and for linear wave equations on time dependent domains \cite{MR4913691,MR4156784}. For a survey of this topic, see \cite{WIDE}.

Motivated by the convexity of \eqref{eq: continuous functional no f} and the convergence of its minimisers to the exact solution of \eqref{eq: semilin wave equation}, we develop and study a higher order, conforming space-time Galerkin method for computing the minimisers of \eqref{eq: continuous functional no f}. As we shall see below, the existence and uniqueness of the discrete approximation can therefore be inferred by that of the continuous problem. 

\subsection{Minimisation for the continuous problem}
\label{sec: minisation}

For a finite time horizon $T>0$, time interval $I:=(0,T]$ and a domain $\Omega\subset\mathbb{R}^d$, $d=1,2,3$, we shall make use of the standard notation for Lebesgue and Sobolev spaces $L^p(\Omega)$ and $W^{k,p}(\Omega)$ with norms $\|\cdot\|_{L^p(\Omega)}$ and $\|\cdot\|_{W^{k,p}(\Omega)}$, (also $|\cdot|_{W^{k,p}(\Omega)}$ denoting the respective seminorm), $k\ge 0$ and $p\in[1,\infty]$, respectively, so that $W^{0,p}(\Omega)=L^p(\Omega)$; when $p=2$ we write $H^k(\Omega)=W^{k,2}(\Omega)$ as is standard, and $H^1_0(\Omega)$ consisting of functions in $H^1(\Omega)$ whose trace on the boundary is zero. Finally, we let $\langle \cdot,\cdot\rangle_{L^2(\Omega)}$ denote the standard $L^2$ inner product. We also consider the Bochner spaces $L^p(I;X)$ with norm
$
    \|u\|_{L^p(I;X)}=\left(\int_0^T||u||_{X}^p\;\text{d}t\right)^{\frac{1}{p}},
$
with $\|\cdot\|_{X}$ denoting the norm of a Banach space $X$. We also define weighted Bochner norms by 
\begin{align*}
    \|u\|_{L^q_{\varepsilon}(I;X)}= \begin{cases}\Big(\int_0^Te^{-\frac{t}{\varepsilon}}\|u\|_{X}^q\;\ud t\Big)^{\frac{1}{q}}, {\qquad q < \infty} \\
    \operatorname*{ess~sup}_{t \in I
    } \;|e^{-\frac{t}{\varepsilon}}\|u(t)\|_{X}|, \qquad  q=\infty.
    \end{cases}
\end{align*}
noting that the finite time horizon $T$ ensures that $\|u\|_{L^q_{\varepsilon}(I;X)}$ is equivalent to $\|u\|_{L^q(I;X)}$. We also define Bochner-Sobolev spaces $H^k(I;X)$ and $H^k_{\varepsilon}(I;X)$ with obvious meaning.  When no confusion is likely to occur, we will abbreviate $\|u\|_{Y(I;W^{k,p}(\Omega))}$ to $\|u\|_{Y(W^{k,p})}$ for $Y\in \{L^q,H^k\}$. For two Banach spaces $Y$ and $Z$, we use the notation $Y\hookrightarrow Z$ to indicate a continuous embedding of $Y$ into $Z$.

Recalling that $D=\Omega\times I$, we define the space
$$\mathcal{U} := H^2(I;L^2(\Omega))\cap L^2(I;H^1_0(\Omega))\cap L^p(D),$$
and we consider the family of functionals $\mathcal{I}_{\varepsilon}:\mathcal{U}\rightarrow \mathbb{R}$ given by
\begin{equation}
    \mathcal{I}_{\varepsilon}(u) = \int_0^Te^{-\frac{t}{\varepsilon}}\Big(\varepsilon^2\|u_{tt}\|_{L^2(\Omega)}^2+\|\nabla u\|_{L^2(\Omega)}^2+\|u\|_{L^p(\Omega)}^p-2\left<f,u\right>_{L^2(\Omega)}\Big)\,\ud t.
    \label{eq: Continuous functional}
\end{equation}

While the homogeneous Dirichlet boundary condition is enforced in the space $\mathcal{U}$, to impose the initial conditions, we minimise $\mathcal{I}_{\varepsilon}$ over the space
\begin{align*}
    \mathcal{K}(u^0,u^1)=\{u \in \mathcal{U}:u(0)=u^0,u_t(0)=u^1\}.
\end{align*}
This set is well defined as $\mathcal{K}(u^0,u^1) \subset H^2(0,T;L^2(\Omega))$ and the point values at $t=0$ exist in the sense of traces. To avoid notational overhead, we set $u^0=0,u^1=0$ and minimise over the space $\mathcal{K}(0,0)$. We stress, however, that the treatment of general initial conditions follows completely analogously. 

For completeness, we prove the existence of a minimiser. 
\begin{theorem} $\mathcal{I}_{\varepsilon}(\cdot)$ has a unique minimiser over $\mathcal{K}(0,0)$.
\end{theorem}
\begin{proof}
$\mathcal{I}_{\varepsilon}(u)$ is strictly convex, with convexity inherited from the convexity of the norms. $\mathcal{K}(0,0)\neq \emptyset$ since $t\mapsto 0 \in \mathcal{K}(0,0)$, and is also convex and closed. To see the latter, let $(v_n)_n\rightarrow v$ with $v_n \in \mathcal{K}(0,0)$. As $\mathcal{K}(0,0)\subset \mathcal{U}$ and $\mathcal{U}$ is a Banach space and hence closed,  $v \in \mathcal{U}\hookrightarrow H^2(L^2)$. Furthermore, since $H^2(L^2)\hookrightarrow C^1(I;L^2(\Omega))$, we have $v(0)=0,v_t(0)=0$. Thus, $v \in \mathcal{K}(0,0)$ and $\mathcal{K}(0,0)$ is closed.

We now show that $\mathcal{I}_{\varepsilon} $ is coercive on $\mathcal{K}(0,0)$. As $u \in \mathcal{K}(0,0)\subset H^2(L^2)$, the fundamental theorem of Calculus in Bochner spaces, (see, e.g., \cite{MR1625845},) along with H\"{o}lder's inequality and Fubini's Theorem give 
\begin{align*}
    \|u(\cdot,t)\|_{L^q(\Omega)}^q&
    \leq\left(\int_0^t(t-s)^{q'}\ud s\right)^{\frac{q}{q'}}\int_0^t||u_{tt}||_{L^q}^q\;\ud s 
    =(q'+1)^{-\frac{q}{q'}}t^{\frac{q}{q'}(q'+1)}\int_0^t||u_{tt}||_{L^q}^q\;\ud s,
\end{align*}
for $q>1$ and $u \in W^{2,q}(L^q)$. Since  $e^{-\frac{t}{\varepsilon}}\int_0^t\|u_{tt}(\cdot,s)\|_{L^q(\Omega)}^q\ud s \leq \int_0^te^{-\frac{s}{\varepsilon}}\|u_{tt}(\cdot,s)\|_{L^q(\Omega)}^q\ud s$, we have
\begin{equation}
\begin{split}
    \int_0^Te^{-\frac{t}{\varepsilon}}\|u\|_{L^q(\Omega)}^q\,\ud t
    &\leq (q'+1)^{-\frac{q}{q'}}\int_0^Tt^{\frac{q}{q'}(q'+1)}\int_0^te^{-\frac{s}{\varepsilon}}||u_{tt}||_{L^q}^q\;\ud s\;\ud t 
    \leq \frac{T^{2q}}{2}\int_0^Te^{-\frac{t}{\varepsilon}}\|u_{tt}\|_{L^q}^q\;\text{d}t.
    \label{eq: second deriv norm Lk}
\end{split}
\end{equation}
In particular, for $q=2$,
\begin{equation}
    \|u\|_{L^2_\varepsilon(I;L^2(\Omega))}^2
    \leq \frac{T^4}{2} \|u_{tt}\|_{L^2_\varepsilon(I;L^2(\Omega))}^2.
\label{eq: second deriv norm}
\end{equation}
Hence the $H^2_{\varepsilon}(L^2)$-seminorm defines a norm on $H^2_{\varepsilon}(L^2)$, with trivial initial data. Standard estimation from below, along with the estimate \eqref{eq: second deriv norm}, gives
\begin{align*}
    \mathcal{I}_{\varepsilon}(u)
    \geq 
    & \  \int_0^Te^{-\frac{t}{\varepsilon}}\Big( \varepsilon^2\|u_{tt}\|_{L^2(\Omega)}^2+\|\nabla u||_{L^2(\Omega)}^2+\|u\|_{L^p(\Omega)}^p\Big)\,\ud t\\
    &-\sqrt{2}T^2\Big(\int_0^Te^{-\frac{t}{\varepsilon}}\|f\|_{L^2(\Omega)}^2\,\ud t\Big)^{\frac{1}{2}}\Big(\int_0^Te^{-\frac{t}{\varepsilon}}\|u_{tt}\|_{L^2(\Omega)}^2\, \ud t\Big)^{\frac{1}{2}}.
\end{align*}
Coercivity then follows from the fact that $\lim_{x\rightarrow \pm\infty} ax^2-bx=+\infty$ for $a,b\geq 0$. Hence, due to the strict convexity of the norm, the functional attains a unique minimiser over the space $\mathcal{K}(0,0)$.
\end{proof}
\begin{remark}
In the infinite time-horizon setting $(T=+\infty)$, including a nonhomogeneous forcing is covered in \cite{MR3902175}, which is a nontrivial extension of the arguments presented in \cite{MR2912711}. 
\end{remark}

\subsection{Euler-Lagrange equation}
Denote the unique minimiser of $\mathcal{I}_{\varepsilon}(\cdot)$ by $u^{\varepsilon} \in \mathcal{K}(0,0)$. The Euler-Lagrange equation associated with $\mathcal{I}_{\varepsilon}$ is given by: find $u^{\varepsilon}\in \mathcal{K}(0,0)$, such that
\begin{equation}
\begin{split}
  A(u^{\varepsilon},v)+B(u^{\varepsilon};v) = f(v), \qquad \text{for all } v \in \mathcal{K}(0,0),
    \label{eq: Continuous Euler-lagrange equation}
\end{split}
\end{equation}
whereby
\begin{align*}
    A(u,v)  :=&\  \int_0^Te^{-\frac{t}{\varepsilon}}\Big(\varepsilon^2\left<u_{tt},v_{tt}\right>_{L^2(\Omega)}+\left<\nabla u,\nabla v\right>_{L^2(\Omega)}\Big)\,\ud t, \\
    B(u;v) :=&\  \frac{p}{2}\int_0^Te^{-\frac{t}{\varepsilon}}\left<|u|^{p-2}u,v\right>_{L^2(\Omega)}\, \ud t, \qquad
    f(v) :=\ \int_0^Te^{-\frac{t}{\varepsilon}}\left<f,v\right>_{L^2(\Omega)}\, \ud t.
\end{align*}
The finite element discretisation below will be based on \eqref{eq: Continuous Euler-lagrange equation}. 

Now, assuming sufficient regularity on the solution $u^{\varepsilon}$, upon integration by parts, we deduce
\begin{equation}
\begin{split}
    0=& \int_0^Te^{-\frac{t}{\varepsilon}}\left<\varepsilon^2u^{\varepsilon}_{tttt}-2\varepsilon u^{\varepsilon}_{ttt}+u^{\varepsilon}_{tt}-\Delta u^{\varepsilon}+\frac{p}{2}|u^{\varepsilon}|^{p-2}u^{\varepsilon}-f,v\right>_{L^2(\Omega)}\!\!\ud t\\
    &+e^{-\frac{T}{\varepsilon}}\left<\varepsilon^2u^{\varepsilon}_{tt}(T),v_t(T)\right>_{L^2(\Omega)} \!\! -e^{-\frac{T}{\varepsilon}}\left<\varepsilon^2u^{\varepsilon}_{ttt}(T)-\varepsilon u^{\varepsilon}_{tt}(T),v(T)\right>_{L^2(\Omega)},
    \label{eq: Almost strong form}
\end{split}
\end{equation}
for all $v \in \mathcal{K}(0,0)$. We conclude formally that $u^{\varepsilon}$ solves the following problem in strong form: 
\begin{equation}
\begin{split}
\varepsilon^2u^{\varepsilon}_{tttt}-2\varepsilon u^{\varepsilon}_{ttt}+u^{\varepsilon}_{tt}-\Delta u^{\varepsilon}+\frac{p}{2}|u^{\varepsilon}|^{p-2}u^{\varepsilon}=f, &\ \text{ in } \Omega \times I, \\
u^{\varepsilon}=0, &\ \text{ on }\partial\Omega \times I, \\
u^{\varepsilon}(\cdot,0)=0,u_t^{\varepsilon}(\cdot,0)=0, &\ \text{ in }\Omega \\
\varepsilon u^{\varepsilon}_{tt}(\cdot,T)=0,\varepsilon^2u_{ttt}^{\varepsilon}(\cdot,T)=0,&\ \text{ in } \Omega.
\end{split}
\label{eq: fourth order time PDE}
\end{equation}
Hence \eqref{eq: fourth order time PDE} is the elliptic regularisation of \eqref{eq: Model equation}. The terminal conditions $\varepsilon u_{tt}^{\varepsilon}(\cdot,T)=\varepsilon^2 u_{ttt}^{\varepsilon}(\cdot,T)=0$ are non-essential. Nevertheless, they deserve special attention as their presence makes this problem non-causal. 
\begin{remark}
This argument can be replicated in the case of nontrivial initial data, $u^0,u^1 \in L^p(\Omega)\cap H^1_0(\Omega)$ and/or for non-homogenenous Dirichlet data. In these cases, minimisation takes place over the respective solution sets instead of $\mathcal{U}$.
\end{remark}
\subsection{Boundary layer and dispersion relation} \label{BL}
Formally, taking $\varepsilon \rightarrow 0^+$ in \eqref{eq: fourth order time PDE} one recovers \eqref{eq: Model equation}, so the behaviour of \eqref{eq: fourth order time PDE} is of interest for $\varepsilon\ll 1$. Applying the scaling $s = \frac{T-t}{\varepsilon}$, we can formally see that the solution $u^\varepsilon$ has a boundary layer of size $\mathcal{O}(\varepsilon)$ at $t=T$. In the layer's inner region the second and third derivatives change rapidly to enforce the additional conditions $u^{\varepsilon}_{tt}(\cdot,T)=0$, $u_{ttt}^{\varepsilon}(\cdot,T)=0$. The outer solution to this perturbation problem is given by \eqref{eq: fourth order time PDE} without the additional terminating conditions, meanwhile the solution inside the boundary layer takes the form:
\begin{equation}
\begin{split}
    u^{\varepsilon}_{ssss}+2u^{\varepsilon}_{sss}+u^{\varepsilon}_{ss}-\varepsilon^2 \Delta u^{\varepsilon} +\frac{p\varepsilon^2}{2}|u^{\varepsilon}|^{p-2}u^{\varepsilon}= \varepsilon^2 f,\\
    u^{\varepsilon}_{ss}(0)=u^{\varepsilon}_{sss}(0)=0.
\end{split}
\label{eq: Inner solution}
\end{equation}
Thus, inside the boundary layer the solution $u^{\varepsilon}$ is a perturbation of an affine function of $s$. 

We seek to quantify the error in approximating the solution of \eqref{eq: Model equation} by \eqref{eq: fourth order time PDE} in terms of $\varepsilon$ as well as the qualitative behaviour of the approximation in terms of the dispersion relation. To that end, we set  $u^{\varepsilon}(\mathbf{x},t) = e^{i(\mathbf{k}\cdot \mathbf{x}-\omega t)}$ into the linear problem \eqref{eq: fourth order time PDE} to deduce the dispersion relation

\begin{align*}
    i \varepsilon \omega^2+\omega \pm |\mathbf{k}| = 0.
\end{align*}
For $\varepsilon\ll1$, writing this as a perturbative solution, we have
\begin{align*}
    \omega = \pm |\mathbf{k}|-i\varepsilon|\mathbf{k}|^2+\mathcal{O}(\varepsilon^2).
\end{align*}
The complex frequency indicates that the solution will decay exponentially with rate $\varepsilon|\mathbf{k}|^2$, since
$
    u^{\varepsilon}(\mathbf{x},t) 
    = e^{-\varepsilon|\mathbf{k}|^2t+\mathcal{O}(\varepsilon^2)}e^{i(\mathbf{k}\cdot \mathbf{x}\pm|\mathbf{k}|t)}.
$
In other words, all sufficiently large frequencies,   $|\mathbf{k}| \gg \varepsilon^{-1/2}$, are damped.  
Comparing this solution to the case $\varepsilon=0$, 
we have
\begin{align*}
    u^{\varepsilon}-u = \left(1-e^{-\varepsilon|\mathrm{k}|^2t + \mathcal{O}(\varepsilon^2)}\right)u,\quad\text{ giving }\quad 
    u^{\varepsilon}-u = \big(\varepsilon|\mathbf{k}|^2t + \mathcal{O}(\varepsilon^2)\big)u,
\end{align*}
to leading order, which indicates (formally) that the error between the solutions of \eqref{eq: Model equation} and \eqref{eq: fourth order time PDE} is of order $\varepsilon$. {blue} We refer to Appendix \ref{sec: Asymptotic analysis} for an analysis involving asymptotic expansions of the boundary layer.

\section{Space-time finite element method}
\label{sec: discrete min}

\begin{figure}[ht]
    \centering
    \begin{subfigure}[b]{0.44\textwidth}
       \centering
    \includegraphics[scale=0.41]{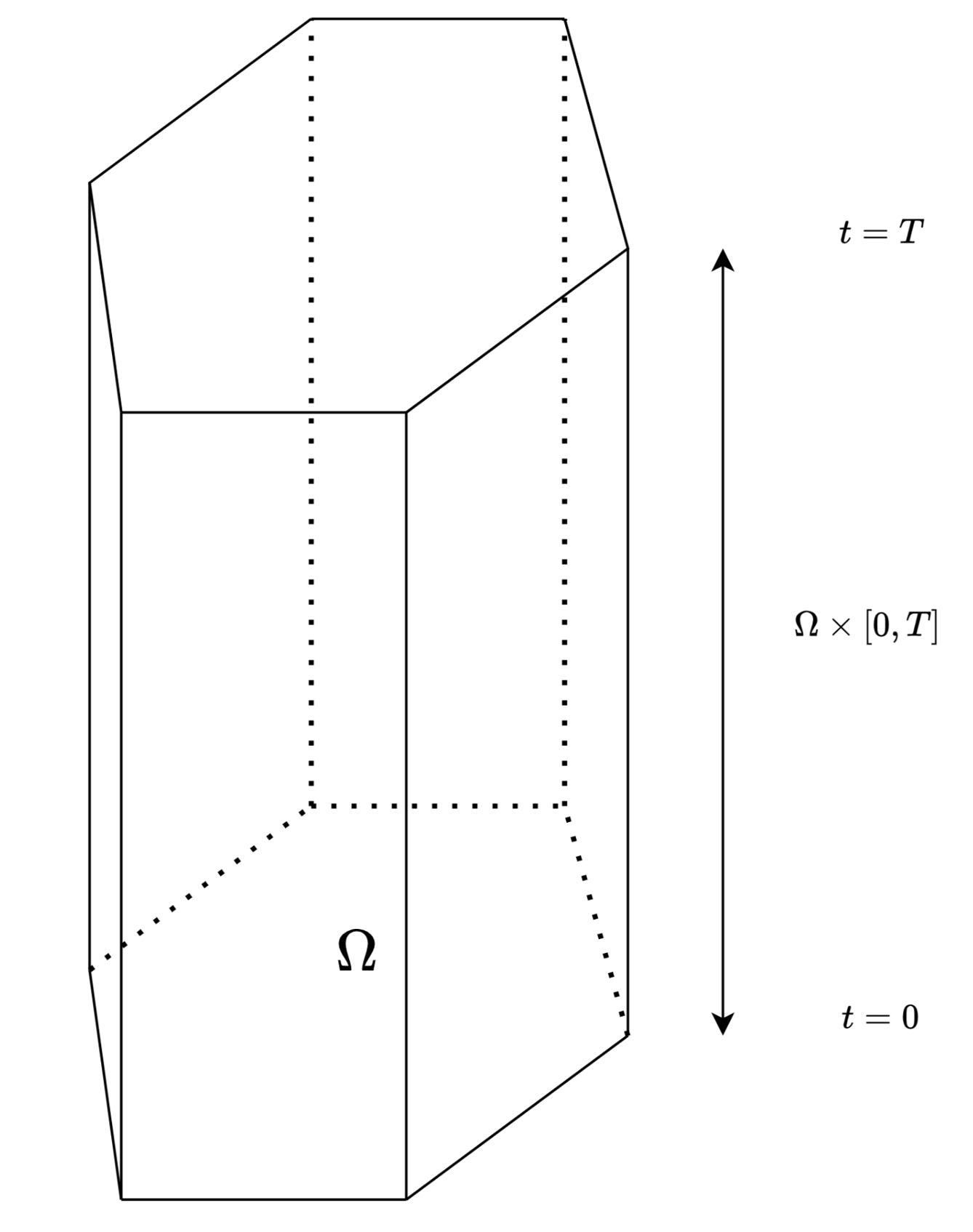}  
        \caption{space-time cylinder $D$}
        \label{fig: space-time cyclinder}
    \end{subfigure}
    \begin{subfigure}[b]{0.44\textwidth}
       \centering
    \includegraphics[scale=0.3]{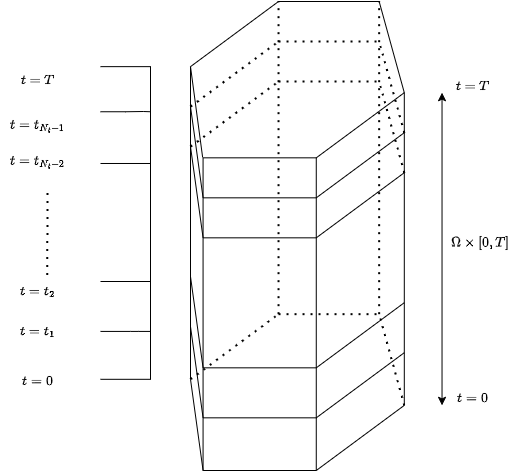}  
        \caption{space-time slabs $\Omega \times I_i$}
        \label{fig: space-time cyclinder with slabs}
    \end{subfigure}
\caption{The space-time cylinder $D$  (\ref{fig: space-time cyclinder}) and its subdivision into space-time slabs (\ref{fig: space-time cyclinder with slabs}). }
\label{fig: space-time cylinder demonstration}
\end{figure}

We subdivide the space-time cylinder $D= \Omega\times I$ into space-time slabs, by dividing the interval $I=(0,T]$ into $N_t$ subintervals $I_i:=(t_i,t_{i+1}]$, $i=0,\cdots,N_t-1$, with $t_0=0,t_{N_t}=T$. We denote $\tau_i = |t_{i+1}-t_i|$ and the timestep as $\tau = \max_{i=0,\cdots,N_{t-1}} \tau_i$ and the collection of these subintervals as $\Delta_{\tau}=\{I_i\}_{i=0}^{N_t-1}$. The space-time slabs $\Omega\times I_{i}$ are illustrated in Figure \ref{fig: space-time cylinder demonstration}.

\begin{figure}[ht]
    \centering
    \begin{subfigure}[b]{0.48\textwidth}
       \centering
    \includegraphics[scale=0.23]{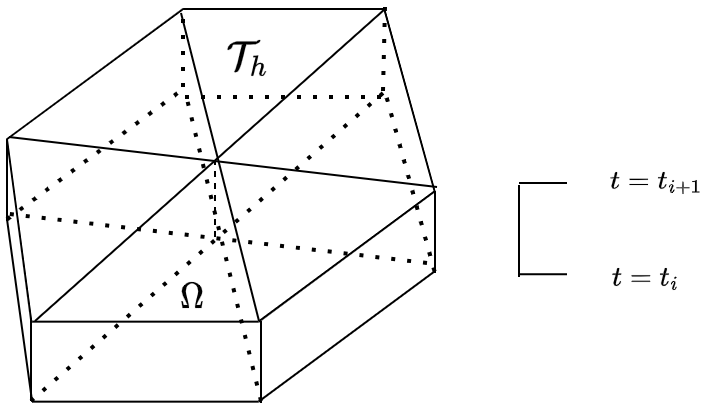}  
        \caption{spatial mesh $\mathcal{T}_h$ on $\Omega \times I_i$ }
        \label{fig: individual slab}
    \end{subfigure}
    \begin{subfigure}[b]{0.49\textwidth}
       \centering
    \includegraphics[scale=0.3]{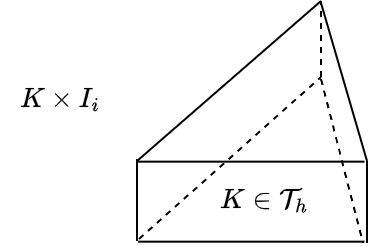}  
        \caption{Prismatic space-time element.}
        \label{fig: individual element}
    \end{subfigure}
\caption{Illustration of the spatial triangulation $\mathcal{T}_{h}$ on the space-time slab $\Omega\times I_i$ in \ref{fig: individual slab} and the individual prismatic elements in \ref{fig: individual element}.}
\label{fig: prismatic elements}
\end{figure}

Next, we place a shape-regular mesh $\mathcal{T}_h = \{K_j\}$ on the polygonal/polyhedral spatial domain $\Omega$, so that $ \Bar{\Omega} = \cup_{j=1}^{N_K}K_j$, $\mathring{K}_n\cap\mathring{K}_m=\emptyset$, $n\neq m$. The space-time domain is, thus, split into discretised space-time slabs as
$    \Bar{\Omega} \times (0,T] = \cup_{i=0}^{N_t-1}\cup_{j=1}^{N_K}K_j\times I_{i}$;
see Figure \ref{fig: prismatic elements}.

Let
\begin{align*}
    V_h &= \{v \in H^1_0(\Omega):v|_K \in \mathcal{P}_s(K), \text{for all } K \in \mathcal{T}_h\},
\end{align*}
with $\mathcal{P}_s(K)$ the space of polynomials of total degree $s$ on $K$. Let $h:=\max_{K\in\mathcal{T}_h}h_K$, $h_K:={\rm diam }(K)$. By construction, $V_h$ is conforming in $H^1_0(\Omega)$.

For temporal discretisation, we use $C^2$-regular splines. For a subdivision $\Delta_{\tau}$ and for $r \in \mathbb{N}$, we define 
\begin{equation}\label{eq: spline defintion}
    \mathcal{S}^{\Delta_{\tau}}(r+2) := \left\{\phi \in C^{2}([0,T]): \phi|_{I_j} \in \mathcal{P}_{r+2}(I_j), j=0,\cdots,N_t-1\right\},
\end{equation}
 so that $\mathrm{dim}(\mathcal{S}^{\Delta_{\tau}}(r+2)) = rN_t+3$. 
 Note that for $r=1$ this corresponds to maximal regularity cubic splines. To enforce the zero initial data, we set $U_{\tau} = \mathcal{S}^{\Delta_{\tau}}(r+2)\cap \{u \in C^1([0,T]): u(0)=0,u_t(0)=0\}$. 

The space-time finite element space is given by $\mathcal{K}_{\tau,h} = U_{\tau}\otimes V_h$. By construction, this is a conforming discretisation since $\mathcal{K}_{\tau,h}\subset\mathcal{K}(0,0)$. 
Therefore, the minimisation problem has a unique solution $u^{\varepsilon}_{\tau,h}$ on $\mathcal{K}_{\tau,h}$, such that
\begin{equation}
\begin{split}
   A(u^{\varepsilon}_{\tau,h},v)+B(u^{\varepsilon}_{\tau,h};v) = f(v), \qquad \text{for all } v \in \mathcal{K}_{\tau,h}.
    \label{eq: Discrete Euler-lagrange equation}
\end{split}
\end{equation}
Thus, by construction the method is well-posed and is also unconditionally stable in approximating $u^\varepsilon$; in particular no CFL condition relating the spatial and temporal meshes is required.

\begin{remark}
\label{rem: Non-triv initial data}
The case of non-trivial initial data can be treated analogously to non-zero essential boundary Dirichlet conditions in finite element analysis. Let us first assume that the initial data is in the discrete space, i.e., $u^0, u^1 \in V_h$. We then look for the solution $u^{\varepsilon}_{\tau,h} \in \tilde{ \mathcal{K}}_{\tau,h} = \mathcal{S}^{\Delta_{\tau}}(r+2) \otimes V_h.$ such that \eqref{eq: Discrete Euler-lagrange equation} holds and $u^{\varepsilon}_{\tau,h}(\cdot,0) = u^0$, $\partial_t u^{\varepsilon}_{\tau,h}(\cdot,0) = u^1$. Again, via the minimisation problem, the unique solution exists. In addition, the error analysis described in the next section proceeds analogously. If the initial data is not in the discrete space, a projection of the initial data onto $V_h$ is required which introduces an additional error. 
\end{remark}
\begin{remark}
The well-posedness does not depend on the tensor product structure of the discrete space. We emphasise that the well-posedness only depends on the subspace property, which would allow for more generic space-time simplicial meshes or even non-polynomial approximation spaces. 
\end{remark}

\section{A priori error analysis}\label{sec: Error rates}

We first establish a quasi-optimality result. In what follows, the dependence of constants on important parameters will be given as arguments.

\begin{theorem}[Quasi-optimality]\label{theorem:ceas lemma}
Let $u^{\varepsilon}$ and $u^{\varepsilon}_{\tau,h}$ solve \eqref{eq: Continuous Euler-lagrange equation} and \eqref{eq: Discrete Euler-lagrange equation}, respectively, and set $\error:=u^{\varepsilon}-u^{\varepsilon}_{\tau,h}$, for brevity. For constants $\alpha\equiv\alpha(p)$, $\beta\equiv\beta(p,C_{\Omega},T,f)$, and $\gamma\equiv\gamma(p)$, we have
\begin{align*}
    &\|\error\|_{E_{\varepsilon}}^2 +\alpha\|\error\|_{L^p_{\varepsilon}(L^p)}^p \\
    &\leq\inf_{v \in \mathcal{K}_{\tau,h}}\left( \|u^{\varepsilon}-v\|_{E_{\varepsilon}}^2 +\beta \|u^{\varepsilon}-v\|_{L^p_{\varepsilon}(L^p)}^{p'}+\gamma\|u^{\varepsilon}-v\|_{L^p_{\varepsilon}(L^p)}^p\right),
\end{align*}
with $p' = \frac{p}{p-1}$ the H\"{o}lder conjugate to $p$, and the energy norm given by
\begin{equation}\label{eq: energy norm}
    \|g\|_{E_{\varepsilon}} :=\sqrt{ A(g,g)} = \Big(\int_0^Te^{-\frac{t}{\varepsilon}}\Big(\varepsilon^2\|g_{tt}\|_{L^2(\Omega)}^2+\|\nabla g\|_{L^2(\Omega)}^2\Big)\,\ud t\Big)^{\frac{1}{2}}.
\end{equation}

\end{theorem}
\begin{proof}
We begin by observing the Galerkin orthogonality property:
\begin{equation}\label{eq: galerkin orthogonality}
    A(\error,v) = -\big(B(u^{\varepsilon};v)-B(u^{\varepsilon}_{\tau,h};v)\big),\quad\text{for all } v \in \mathcal{K}_{\tau,h}.
\end{equation}
For an arbitrary $v \in \mathcal{K}_{\tau,h}$, we set
$
    \eta := u^{\varepsilon} -v$ and  
$    \zeta := v-u^{\varepsilon}_{\tau,h}$,
so that $\error  = \eta + \zeta $, and
\begin{equation}
\begin{aligned}
\|\error\|_{E_{\varepsilon}}^2 
      =&\  A(\error,\eta)-(B(u^{\varepsilon};\zeta)-B(u^{\varepsilon}_{\tau,h};\zeta))\\
     = &\ A(\error,\eta)-\big(B(u^{\varepsilon};\error)-B(u^{\varepsilon}_{\tau,h};\error)\big)+\big(B(u^{\varepsilon};\eta)-B(u^{\varepsilon}_{\tau,h};\eta)\big)\\
\equiv&\ A(\error,\eta)- (I)+(II),
\end{aligned}\label{eq: nonlinear error equation}
\end{equation}
employing Galerkin orthogonality \eqref{eq: galerkin orthogonality}.
We bound each term separately.

We employ the monotonicity of the map $x \mapsto |x|^{p-2}x$ to conclude that there exists a $C_p>0$, such that
\begin{align}
\label{eq: monotonicity}
  (I)=&\ \frac{p}{2}\int_0^T e^{-\frac{t}{\varepsilon}} \left<|u^{\varepsilon}|^{p-2}u^{\varepsilon}-|u_{\tau,h}^{\varepsilon}|^{p-2}u_{\tau,h}^{\varepsilon},\error\right>_{L^2(\Omega)}\,\ud t\\
  \geq&\  \frac{pC_p}{2}\int_0^T e^{-\frac{t}{\varepsilon}} \|\error\|_{L^p(\Omega)}^p\,\ud t =\frac{pC_p}{2}\|\error\|_{L^p_{\varepsilon}(L^p)}^p.
\end{align}

Now, using a local Lipschitz estimate of the map $x \mapsto |x|^{p-2}x$, we ascertain that there exists a constant $C_1>0$, such that
\begin{align*}
    \left\lvert|u^{\varepsilon}|^{p-2}u^{\varepsilon}-|u^{\varepsilon}_{\tau,h}|^{p-2}u^{\varepsilon}_{\tau,h}\right\rvert \leq C_1(|u^{\varepsilon}|+|u_{\tau,h}^{\varepsilon}|)^{p-2}|\error|,
\end{align*}
so that
\begin{align*}
   |(II)|
    \leq &\ \frac{p C_1}{2}\int_0^T e^{-\frac{t}{\varepsilon}}\|\;|u^{\varepsilon}|+|u_{\tau,h}^{\varepsilon}|\;\|_{L^p(\Omega)}^{p-2}\|\error\|_{L^p(\Omega)}||\eta||_{L^p(\Omega)}\,\ud t\\
    \leq&\ \frac{pC_1}{2}||\;|u^{\varepsilon}|+|u^{\varepsilon}_{\tau,h}|\;\|^{p-2}_{L^p_{\varepsilon}(L^p)}\|\error\|_{L^p_{\varepsilon}(L^p)}\|\eta\|_{L^p_{\varepsilon}(L^p)},\label{eq:Lp errors}
\end{align*}
where we have applied the generalised H\"{o}lder's inequality with exponents $p,p, \frac{p}{p-2}$ for the first inequality and the elementary identity $e^{-\frac{t}{\varepsilon}} = e^{-\frac{t}{p\varepsilon}}\cdot e^{-\frac{t}{p\varepsilon}}\cdot e^{-\frac{t}{\varepsilon}\left(\frac{p-2}{p}\right)}$, along with H\"{o}lder's inequality again for the second inequality.

Further, Young's inequality gives $|A(\eta+\zeta,\eta)|\leq \frac{1}{2}\|\error\|_{E_{\varepsilon}}^2+\frac{1}{2}\|\eta\|_{E_{\varepsilon}}^2$. For the nonlinear term, we work as follows:
\begin{align*}
    \|\,|u^{\varepsilon}|+|u^{\varepsilon}_{\tau,h}|\,\|^{p-2}_{L^p_{\varepsilon}(L^p)} &\leq  \|2|u^{\varepsilon}|+|\error|\;\|^{p-2}_{L^p_{\varepsilon}(L^p)}
    \leq 2^{2p-5}\big(\|u^{\varepsilon}\|_{L^p_{\varepsilon}(L^p)}^{p-2}+ \|\error\|_{L^p_{\varepsilon}(L^p)}^{p-2}\big).
\end{align*}
Using the above bounds on \eqref{eq: nonlinear error equation}, results in
\begin{equation}\label{eq: nonlinear bound 3}
\begin{split}
&\frac{1}{2}\|\error\|_{E_{\varepsilon}}^2+\frac{pC_p}{2}\|\error\|_{L^p_{\varepsilon}(L^p)}^p \\
\leq&\ \frac{1}{2}\|\eta\|_{E_{\varepsilon}}^2
+2^{2p-6}pC_1 \big(\|u^{\varepsilon}\|_{L^p_{\varepsilon}(L^p)}^{p-2}+ \|\error\|_{L^p_{\varepsilon}(L^p)}^{p-2}\big)\|\error\|_{L^p_{\varepsilon}(L^p)}\|\eta\|_{L^p_{\varepsilon}(L^p)}.
\end{split}
\end{equation}

Now, using Young's inequality with parameters $\sigma_1,\sigma_2$ and exponents $p,p' = \frac{p}{p-1}$, we obtain
\begin{equation}
\begin{split}
&\frac{1}{2}\|\error\|_{E_{\varepsilon}}^2+\Big(\frac{pC_p}{2}-C_12^{2p-6}(\sigma_1^p+\sigma_2^p)\Big)\|\error\|_{L^p_{\varepsilon}(L^p)}^p\\
\leq &\ \frac{1}{2}\|\eta\|_{E_{\varepsilon}}^2+\frac{pC_12^{2p-6}}{p'\sigma_1^{p'}}\|u^{\varepsilon}\|_{L^p_{\varepsilon}(L^p)}^{p'(p-2)}\|\eta\|_{L^p_{\varepsilon}(L^p)}^{p'} +\frac{pC_12^{2p-6}}{p'\sigma_2^{p'}}\|\eta\|_{L^p_{\varepsilon}(L^p)}^p.
\label{eq: nonlinear bound 4}
\end{split}
\end{equation}
Next choose $\sigma_1,\sigma_2$ independent of $\varepsilon$ small enough, such that
\begin{align*}
  \frac{pC_p}{4}\ge C_12^{2p-6}(\sigma_1^p+\sigma_2^p).
\end{align*}
With this choice of $\sigma_1,\sigma_2$ define the constants
\begin{align*}
    \alpha = \frac{pC_p}{2},\quad  \rho = \frac{pC_12^{2p-5}}{p'\sigma_1^{p'}}\|u^{\varepsilon}\|_{L^p_{\varepsilon}(L^p)}^{p'(p-2)},\quad  \gamma=\frac{pC_12^{2p-5}}{p'\sigma_2^{p'}},
\end{align*}
thereby concluding from \eqref{eq: nonlinear bound 4}, the estimate 
\begin{equation}
\|\error\|_{E_{\varepsilon}}^2+\alpha\|\error\|_{L^p_{\varepsilon}(L^p)}^p\\
\leq  \|\eta\|_{E_{\varepsilon}}^2+\rho \|\eta\|_{L^p_{\varepsilon}(L^p)}^{p'} +\gamma\|\eta\|_{L^p_{\varepsilon}(L^p)}^p.
\end{equation}
{To complete the proof, we show that the constant $\rho$ depends only on the data, by establishing that $u^{\varepsilon}$ is bounded in $L^{p}_{\varepsilon}(I;L^p(\Omega))$ by $\|f\|_{L^2(L^2)}$. To that end, setting $v=u^{\varepsilon}$ in \eqref{eq: Continuous Euler-lagrange equation} and employing Cauchy-Schwarz and Young's inequalities, we deduce
\[
\begin{aligned}
\|\nabla u\|_{L^2_{\varepsilon}(L^2)}^2+ \|u^{\varepsilon}\|_{L^p_{\varepsilon}(L^p)}^p\leq&\ A( u^{\varepsilon},u^{\varepsilon})+B(u^{\varepsilon};u^{\varepsilon})
 \leq \frac{C_{\Omega}^2}{2}\|f\|^2_{L^2_{\varepsilon}(L^2)}+\frac{1}{2}\int_0^Te^{-\frac{t}{\varepsilon}}\|\nabla u^{\varepsilon}\|_{L^2(\Omega)}^2\;\text{d}t,
\end{aligned}
\]
 where $C_{\Omega}$ is the Poincaré inequality constant, giving
 \begin{equation}
    \|u^{\varepsilon}\|_{L^p_{\varepsilon}(L^p)}^p\leq \frac{C_{\Omega}^2}{2}\|f\|_{L^2_{\varepsilon}(L^2)}^2 \leq \frac{C_{\Omega}^2}{2}\|f\|_{L^2(L^2)}^2,
    \label{eq: f bound on u}
\end{equation}
and hence 
 $   \rho\le \beta :=  \frac{pC_12^{2p-6}}{p'\sigma_1^{p'}}\Big(\frac{C_{\Omega}^2}{2}\|f\|_{L^2(L^2)}^2\Big)^{\frac{p-2}{p-1}}.$
}
\end{proof}
\begin{remark}
If $p=2$, so that $p'=p$, \eqref{eq: Model equation} corresponds to a linear Klein-Gordon equation. Then, Theorem \ref{theorem:ceas lemma} holds with $\alpha=1$, $\beta+\gamma =1 $, upon simply modifying the energy norm to $\|f\|_{KG}^2 =\|f\|_{E_{\varepsilon}}^2+\|f\|_{L^2_{\varepsilon}}^2$, to obtain $\|e_{\varepsilon}\|_{KG}= \inf_{v \in \mathcal{K}_{\tau,h}}\|u^{\varepsilon}-v\|_{KG}$. Moreover, for the linear wave equation, corresponding to $p=0$, the bound in Theorem \ref{theorem:ceas lemma} reduces to $\|e_{\varepsilon}\|_{E_{\varepsilon}}= \inf_{v\in\mathcal{K}_{\tau,h}}\|u^{\varepsilon}-v\|_{E_{\varepsilon}}$.
\end{remark}
\begin{remark}
A power $p$ in the statement of Theorem \ref{theorem:ceas lemma} may be recovered through the use of quasi-norms as defined in \cite{MR1276708}. In this setting, one may obtain that 
\begin{align*}
    \|e_{\varepsilon}\|_{E_{\varepsilon}}^2+\|e_{\varepsilon}\|_{L^2_{\varepsilon}((u,p))}^2 \leq \Gamma\inf_{v\in \mathcal{K}_{\tau,h}}(\|u^{\varepsilon}-v\|_{E_{\varepsilon}}^2+\|u^{\varepsilon}-v\|_{L^2_{\varepsilon}((u,p))}^2),
\end{align*}
for a constant $\Gamma$ independent of $\varepsilon$. Here the quasi-norm is defined as $\|f\|_{L^2_{\varepsilon}((w,p))}^2 =\int_0^Te^{-\frac{t}{\varepsilon}}\|f\|_{(w,p)}^2\;\text{d}t$, with $\|f\|_{(w,p)}^2 =\int_{\Omega}|f|^2(|w|+|f|)^{p-2}\;\text{d}x$. The proof of this statement follows analogously from the proof of \cite[Theorem 4.6]{MR4780834}, which uses the boundedness and monotonicity properties of the nonlinearity in the context of the quasi-norms; this is omitted for brevity.
\end{remark}

\begin{theorem}[Quasi-optimality for `small' nonlinearities]\label{cor: small nonlin cea} {With $p>2$,} suppose that either $d=1,2$ or if $d \geq 3$, that $p \leq 2+\frac{4}{d}$. Then, there exists a constant $G(p,T,\varepsilon,f)$, with $G = \mathcal{O}(e^{\frac{T}{\varepsilon}})$, such that 
\begin{equation}
\begin{split}
    \|\error\|_{E_{\varepsilon}}\leq \min&\left\{\left(1+G\right)\inf_{v \in \mathcal{K}_{\tau,h}}\|u^{\varepsilon}-v\|_{E_{\varepsilon}},\right. \\
    &\left.\inf_{v \in \mathcal{K}_{\tau,h}}\left( \|u^{\varepsilon}-v\|_{E_{\varepsilon}}^2 +\beta \|u^{\varepsilon}-v\|_{L^p_{\varepsilon}(L^p)}^{p'}+\gamma\|u^{\varepsilon}-v\|_{L^p_{\varepsilon}(L^p)}^p\right)^{\frac{1}{2}}\right\},
\end{split}
\end{equation}
\label{col: small nonlinearity}
with $\beta,\gamma$ independent of $\varepsilon$ as given in Theorem \ref{theorem:ceas lemma}.
\end{theorem}
The second constituent inequality in \eqref{col: small nonlinearity} follows directly from Theorem \ref{theorem:ceas lemma}; the proof of the first inequality, requires a technical result on interpolation in the $\varepsilon$-weighted spaces.

\begin{lemma}\label{lemma: interpolation estimates}
    Suppose that either $d=1,2$ and $2\leq p \leq \infty$, or $d \geq 3$ and $2\le p \leq 2+\frac{4}{d}$. Then, there exists a constant $C_2(\varepsilon,T,\Omega,p,d)>0$, such that  $\|w\|_{L^p_{\varepsilon}(L^p)}\leq C_2(\varepsilon,T,\Omega,p,d)\|w\|_{E_{\varepsilon}}$, for any $w \in \mathcal{K}(0,0)$.
\end{lemma}
\begin{proof}
First, we establish two embedding estimates on the $\varepsilon-$weighted spaces. For $p < m$ and $w \in L^p_{\varepsilon}(L^p(\Omega))\cap L^m_{\varepsilon}(L^m(\Omega))$, we show
\begin{equation}
    \|w\|_{L^p_{\varepsilon}(L^p)}\leq \left(\varepsilon|\Omega|\right)^{\frac{1}{p}-\frac{1}{m}}\|w\|_{L^m_{\varepsilon}(L^m)}.
    \label{eq: weighted holder estimate}
\end{equation}
To begin, for each $t$, we have $\|u\|_{L^p(\Omega)}(t)\leq |\Omega|^{\frac{1}{p}-\frac{1}{m}}\|u\|_{L^m(\Omega)}(t)$, by H\"{o}lder's inequality. Hence,
\begin{align*}
    \|u\|_{L^p_{\varepsilon}(L^p)}^p &\leq |\Omega|^{1-\frac{p}{m}}\int_0^Te^{-\frac{t}{\varepsilon}}\|u\|_{L^m_{\varepsilon}}^p\;\text{d}t
    \leq|\Omega|^{1-\frac{p}{m}}\Big(\int_0^Te^{-\frac{t}{\varepsilon}}\;\text{d}t\Big)^{1-\frac{p}{m}}\|u\|_{L^m_{\varepsilon}(L^m)}^p,
\end{align*}
by H\"{o}lder's inequality in time with exponents $\frac{m}{p}>1$ and $\frac{m}{m-p}$. Noting that the integral is bounded above by $\varepsilon$ and taking the $p^{th}$ root gives \eqref{eq: weighted holder estimate}. 

For the second estimate, for $w \in \mathcal{K}(0,0)$, we show that
\begin{equation}
    \|w\|_{L^{\infty}_{\varepsilon}(L^2)} \leq \sqrt{2\varepsilon}\|w\|_{E_{\varepsilon}},
    \label{eq: L^infty bound on E_{epsilon}}
\end{equation}
with $\|\cdot\|_{E_{\varepsilon}}$ as in \eqref{eq: energy norm}.
As $w \in \mathcal{K}(0,0)$, $w_{tt}\in L^2_{\varepsilon}(L^2)$ and it has trivial initial data. By the the fundamental theorem of Calculus in Bochner spaces, (see, e.g., \cite{MR1625845}), we may write $w(x,t) = \int_0^t(t-s)w_{tt}(x,s)\;\text{d}s$. Using trivial estimation from above, we have
\begin{align*}
    e^{-\frac{t}{\varepsilon}}\|w(t)\|_{L^2}
    &\leq \left(e^{-\frac{2t}{\varepsilon}}\int_0^t\;e^{\frac{s}{\varepsilon}}(t-s)^2\text{d}s\right)^{\frac{1}{2}}\|w_{tt}\|_{L^2_{\varepsilon}(0,t;L^2)}   
    \leq  \sqrt{2}\varepsilon^{\frac{3}{2}}\|w_{tt}\|_{L^2_{\varepsilon}(0,t;L^2)},
\end{align*}
Hence,
 $    \|w\|_{L^{\infty}_{\varepsilon}(L^2)}\leq \sqrt{2\varepsilon}\big(\varepsilon\|w_{tt}\|_{L^2_{\varepsilon}(0,T;L^2)}\big) 
   \leq \sqrt{2\varepsilon}\|w\|_{E_{\varepsilon}}$.

With these two estimates in hand, we establish an interpolation estimate on the weighted $L^m_{\varepsilon}(L^m)$ spaces. Observe that $E_{\varepsilon}\hookrightarrow L^2_{\varepsilon}(H^1)\hookrightarrow L^2_{\varepsilon}(L^{2^{*}})$, for $2^*=\frac{2d}{d-2}$, by the Sobolev embedding. We aim to interpolate $L^m_{\varepsilon}(L^m)$ between $L^2_{\varepsilon}(L^{2^*})$ and $L^{\infty}_{\varepsilon}(L^2)$, as $E_{\varepsilon}$ will control both by \eqref{eq: L^infty bound on E_{epsilon}}. Let $2\leq m \leq 2^{*}$ and $\theta \in (0,1)$, both to be determined later.  Applying a standard $L^m$ interpolation estimate in the spatial variable \cite[Theorem~4.6]{MR2759829}, for each $t$ we have $\|w\|_{L^m(\Omega)}\leq \|w\|_{L^2(\Omega)}^{\theta}\|w\|_{L^{2^*}(\Omega)}^{1-\theta}$, where
\begin{align}
    \frac{1}{m} = \frac{\theta}{2}+\frac{1-\theta}{2^*}.
    \label{eq: interpolation est}
\end{align}
Multiplying by the weight $e^{-\frac{t}{\varepsilon}}$ and integrating, we obtain
\begin{align*}
    \|w\|_{L^m_{\varepsilon}(L^m)}^m &\leq \int_0^Te^{-\frac{t}{\varepsilon}}\|w\|_{L^2}^{m\theta}\|w\|_{L^{2^{*}}}^{m(1-\theta)}\;\text{d}t
    \leq \|w\|_{L^{\infty}_{\varepsilon}(L^2)}^{m\theta}\int_0^Te^{-\frac{t}{\varepsilon}(1-\theta m)}\|u\|_{L^{2^{*}}}^{m(1-\theta)}\;\text{d}t \\
    &\leq e^{\frac{\theta m T}{\varepsilon}}\|w\|_{L^{\infty}_{\varepsilon}(L^2)}^{m\theta} \|u\|_{L^{m(1-\theta)}_{\varepsilon}(L^{2^{*}})}^{m(1-\theta)}.
\end{align*}
 We now choose $m(1-\theta)=2$. With this choice, solving for $\theta, m$ in \eqref{eq: interpolation est} enforces $\theta = 1-\frac{2}{m}$ and $m = 2+\frac{4}{d}$.
Hence,
\begin{equation}
    \|w\|_{L^m_{\varepsilon}(L^m)}\leq e^{\frac{T}{\varepsilon}(1-\frac{2}{m})}\|w\|_{L^{\infty}_{\varepsilon}(L^2)}^{1-\frac{2}{m}}\|w\|_{L^{2}_{\varepsilon}(L^{2^{*}})}^{\frac{2}{m}}.
    \label{eq: interpolation ans}
\end{equation}
Combining these estimates for $p \leq m=2+\frac{4}{d}$, we conclude
\begin{align*}
    \|w\|_{L^p_{\varepsilon}(L^p)}&\leq(\varepsilon|\Omega|)^{\frac{1}{p}-\frac{1}{m}}\|w\|_{L^m_{\varepsilon}(L^m)}, && \text{(from \eqref{eq: weighted holder estimate}})\\
    &\leq(\varepsilon|\Omega|)^{\frac{1}{p}-\frac{1}{m}} e^{\frac{T}{\varepsilon}{\left(1-\frac{2}{m}\right)}}\|w\|_{L^{\infty}_{\varepsilon}(L^2)}^{1-\frac{2}{m}}\|w\|_{L^{2}_{\varepsilon}(L^{2^*})}^{\frac{2}{m}}, && (\text{from } \eqref{eq: interpolation ans})\\
    &\leq (\varepsilon|\Omega|)^{\frac{1}{p}-\frac{1}{m}}e^{\frac{T}{\varepsilon}(1-\frac{2}{m})}(\sqrt{2\varepsilon})^{1-\frac{2}{m}}C^{\frac{2}{m}}\|w\|_{E_{\varepsilon}}, &&(\text{from } \eqref{eq: L^infty bound on E_{epsilon}}) \\
    &\equiv C_2(\varepsilon,T,\Omega,p,d)\|w\|_{E_{\varepsilon}},
\end{align*}
where $C\equiv C(d,m,\Omega)$ is the constant corresponding to the Sobolev embedding $H^1({\Omega})\hookrightarrow L^{2^{*}}(\Omega)$.
\end{proof}

\begin{proof}[Proof of Theorem \ref{col: small nonlinearity}]
We first establish that $u^{\varepsilon}_{\tau,h}$ is bounded in $L^{p}_{\varepsilon}(L^p)$ independently of $\tau,h$. Setting $v = u^{\varepsilon}_{\tau,h}$ in \eqref{eq: Discrete Euler-lagrange equation}, by the same argument, we have that $u^{\varepsilon}_{\tau,h}$ satisfies the same bound as \eqref{eq: f bound on u}. In particular, one has
\begin{equation}
    \|u^{\varepsilon}\|_{L^p_{\varepsilon}(L^p)}^{p-2}+\|u_{\tau,h}^{\varepsilon}\|_{L^p_{\varepsilon}(L^p)}^{p-2}\leq C_3(T,C_{\Omega},p,\|f\|_{L^2(L^2)}).
    \label{eq: constant 4}
\end{equation}
Starting from \eqref{eq: nonlinear bound 3}, using Lemma \ref{lemma: interpolation estimates}, (since $u^{\varepsilon},u^{\varepsilon}_{\tau,h},v \in \mathcal{K}(0,0)$) we may replace the $L^p$ space-time norms with the energy norm, assuming that $p \leq 2+\frac{4}{d}$ and apply H\"{o}lder's inequality on the linear term instead of Young's inequality, to deduce
\begin{equation}\label{eq: nonlinear bound 6}
\begin{split}
    \|\error\|_{E_{\varepsilon}}^2 \leq&\  \|\error\|_{E_{\varepsilon}}\|\eta\|_{E_{\varepsilon}}+ 2^{p-4}C_1C_2^2 p\left(\|u^{\varepsilon}\|_{L^p_{\varepsilon}(L^p)}^{p-2}+\|u^{\varepsilon}_{\tau,h}\|_{L^p_{\varepsilon}(L^p)}^{p-2}\right)\|\eta\|_{E_{\varepsilon}}\|\error\|_{E_{\varepsilon}},
\end{split}
\end{equation}
where we have also dropped the $L^p$ norm on the left hand side.
Applying \eqref{eq: constant 4}, and setting
$    G \equiv G(p,T,\varepsilon,f) = 2^{p-4}C_1C_2^2C_3p$,
the result follows.
\end{proof}
\begin{remark}
    Note that $C_2\sim e^{\frac{T}{\varepsilon}}$ indicating that the constant $G$ resulting from the proof becomes \emph{exponentially large} as $\varepsilon\to 0$. It is not clear if this is the best possible constant achievable.  
\end{remark}
Note that Theorems \ref{theorem:ceas lemma} and \ref{cor: small nonlin cea} pose no restriction on the choice of spatial approximation space, as long as it is conforming. Having established quasi-optimality, our goal is to now prove convergence rates for the method \eqref{eq: Discrete Euler-lagrange equation}. To that end, we establish best approximation results on $\mathcal{K}_{\tau,h}$ in the weighted norms.

For the temporal discretisation, given a Banach space $X$, we define 
$\mathcal{J_{\tau}}:H^2(0,T;X)\rightarrow U_{\tau}\otimes X$, by
\begin{equation}
    \mathcal{J_{\tau}}[u](t) := \int_0^t\int_0^s\mathcal{I}_{\tau}^r[u_{tt}]\;\text{d}w\;\text{d}s,
    \label{eq: temporal interpolant def}
\end{equation}
where $\mathcal{I}_{\tau}^r[v]$ is the piecewise polynomial Lagrange interpolation operator of degree $r$ with respect to the time variable that takes values in $X$,

\begin{lemma}[Interpolation on $\varepsilon$-weighted spaces]
Let $v\in W^{k,p}_{\varepsilon}(X)$, $k\in\mathbb{N}$ and $p\ge 1$. For $r+1\geq k$, define the function $D_p:\mathbb{R}_+\to \mathbb{R}_+$, by
\begin{equation*}
    D_{p}(z) = (p-1)^{\frac{p-1}{p}}z^{-1}\left(e^{\frac{z}{p-1}}-1\right)^{\frac{p-1}{p}}\left(1-e^{-z}\right)^{\frac{1}{p}}.
\end{equation*} 
Then, the following bound holds:
\begin{equation}
    \|v-\mathcal{I}_{\tau}^r[v]\|_{L^p_{\varepsilon}(X)}\leq \tau^{k}C_{r,k}D_p\left(\frac{\tau}{\varepsilon}\right)|v|_{W^{k,p}_{\varepsilon}(X)}.
\label{eq: time Lp error}
\end{equation}
\end{lemma}
\begin{proof}
By definition,
\begin{align*}
    \|v-\mathcal{I}_{\tau}^r[v]\|_{L^p_{\varepsilon}(X)}^p 
    &= \sum_{i=0}^{N_t-1}\int_{t_i}^{t_{i+1}}e^{-\frac{t}{\varepsilon
    }}\|v-\mathcal{I}_{\tau}^r[v]\|_X^p\;\text{d}t 
    \equiv \sum_{i=0}^{N_t-1}E_i.
\end{align*}
For each $E_i$, the interpolant can be expressed in terms of the Lagrange basis functions: let $t_{i} = \xi_{i,0}<\xi_{i,1}<\cdots<\xi_{i,r}=t_{i+1}$ be the nodes associated with the Lagrange basis functions $\{L_{i,j}\}_{j=0,\cdots,r}$ in the subinterval $[t_{i},t_{i+1}]$, giving
\begin{align*}
    \mathcal{I}_{\tau}^r[v](t) = \sum_{j=0}^r v(\xi_{i,j})L_{i,j}(t),
\end{align*}
for $t\in [t_{i},t_{i+1}]$.
Since the Lagrange interpolant is exact for polynomials of degree $k-1\leq r$, the Peano Kernel Theorem (see, e.g., \cite[\S~13.2]{MR2796928}) implies
\begin{align*}
    v(t)-\mathcal{I}_{\tau}^r[v](t) &= \frac{1}{(k-1)!}\int_{t_{i}}^{t_{i+1}}\Big((t-s)_{+}^{k-1}-\sum_{j=0}^rL_{i,j}(t)(\xi_{i,j}-s)_{+}^{k-1}\Big)v^{(k)}(s)\;\text{d}s, 
\end{align*}
with $(\cdot)^{k-1}_{+} := \max\{(\cdot)^{k-1},0\}$. Thus,
\begin{align*}
    \|v(t)-\mathcal{I}_{\tau}^r[v](t)\|_{X}^p \leq 
    & \ 
    \bigg(\frac{1}{(k-1)!}\int_{t_i}^{t_{i+1}}e^{\frac{s}{\varepsilon p}}\Big|(t-s)_{+}^{k-1}-\sum_{j=0}^rL_{i,j}(t)(\xi_{i,j}-s)_{+}^{k-1}\Big|e^{-\frac{s}{\varepsilon p}}\|v^{(k)}(s)\|_{X}\;\text{d}s\bigg)^p\\
    \leq&\ \bigg(\int_{t_i}^{t_{i+1}}e^{\frac{sp'}{\varepsilon p}}\Big|(t-s)_{+}^{k-1}-\sum_{j=0}^rL_{i,j}(t)(\xi_{i,j}-s)_{+}^{k-1}\Big|^{p'}\;\text{d}s\bigg)^{\frac{p}{p'}}\frac{\|v^{(k)}\|_{L^p_{\varepsilon}([t_{i},t_{i+1}];X)}^p}{(k-1)!^p} \\
    \leq& \bigg(\frac{\left(1+\Lambda_{r,i}\right)\tau_i^{(k-1)}}{(k-1)!}\bigg)^p\left(\int_{t_i}^{t_{i+1}}e^{\frac{s}{\varepsilon}\frac{p'}{p}}\;\text{d}s\right)^{\frac{p}{p'}}\|v^{(k)}\|_{L^p_{\varepsilon}([t_{i},t_{i+1}];X)}^p,
\end{align*}
with $\Lambda_{r,i}$ the Lebesgue constant associated with the Lagrange basis on $[t_i,t_{i+1}]$. Multiplying by $e^{-\frac{t}{\varepsilon}}$ and integrating yields
\begin{align*}
   E_i\le  \bigg(\frac{\left(1+\Lambda_{r,i}\right)\tau_i^{k}}{(k-1)!}D_{p}\left(\frac{{\tau_i}}{\varepsilon}\right)\bigg)^p\|v^{(k)}\|_{L^p_{\varepsilon}([t_{i},t_{i+1}];X)}^p,
\end{align*}
upon noting that
\begin{align*}
    \left(\int_{t_i}^{t_{i+1}}e^{\frac{s}{\varepsilon}\frac{p'}{p}}\;\text{d}s\right)^{\frac{p}{p'}}\int_{t_i}^{t_{i+1}}e^{-\frac{t}{\varepsilon}}\;\text{d}t
    &= \tau_i^p D_{p}^p\left(\frac{{\tau_i}}{\varepsilon}\right).
\end{align*}
Note that $D_p(\cdot)\geq 1$, it is strictly increasing and $\lim_{x\rightarrow \infty}D_{p}(x) = +\infty$. Also, the growth of $D_p$ as $\frac{\tau_i}{\varepsilon}\rightarrow \infty$ is exponential. 

Hence with $\Lambda_{r}:= \max_{i=0,\cdots,N_{t-1}}\{\Lambda_{r,i}\}$, the above yield
\begin{align*}
    \|v-\mathcal{I}_{\tau}^r[v]\|_{L^p_{\varepsilon}(X)} &\leq \tau^{k}\frac{\left(1+\Lambda_{r}\right)}{(k-1)!}D_p\left(\frac{\tau}{\varepsilon}\right)|v|_{W^{k,p}_{\varepsilon}(X)},
\end{align*}
from where the result follows.
\end{proof}

In space we assume the existence of a projection operator $\Pi_h\colon C^0(\overline{\Omega}) \to V_h$ such that for some constant $c> 0$
\begin{equation}    
\label{eq: SZ error}
\begin{aligned}
\|v-\Pi_h[v]\|_{H^1(\Omega)}&\leq ch^{l}|v|_{H^{l+1}(\Omega)}, &  0&\leq l \leq s, \quad v \in  H^{l+1}(\Omega)\\
\|v-\Pi_h[v]\|_{L^p(\Omega)}&\leq ch^{l+1}|v|_{W^{p,l+1}(\Omega)},&  -1&\leq l \leq s, \quad v \in W^{p,l+1}(\Omega).
\end{aligned}
\end{equation}
For example, for Lagrange elements, $\Pi_h$ can be chosen as the Lagrange interpolation operator; see \cite[(4.4.20) Theorem]{MR2373954}. The second inequality with $l = -1$  also implies the stability of $\Pi_h$  in $L^p$
\begin{equation}
  \label{SZ_stab} 
  \|\Pi_h w\|_{L^{p}(\Omega)}\le 2c  \|w\|_{L^{p}(\Omega)}, \qquad \text{for all }w\in C^0(\overline{\Omega}).
\end{equation}

The space-time projection operator $\mathcal{J}_{\tau,h}: \mathcal{K}(0,0)\rightarrow \mathcal{K}_{\tau,h}$ is then defined as $\mathcal{J}_{\tau,h}:=\mathcal{J}_{\tau}\circ \Pi_h$.

\begin{lemma}[Space-time projection error]\label{lemma: projection error} Suppose $w \in \mathcal{K}(0,0)$ has additional regularity such that $w\in H^{k+2}_{\varepsilon}(I;H^1(\Omega))\cap W^{k+2,p}_{\varepsilon}(I;L^p(\Omega))$ and 
$w\in H^2_{\varepsilon}(I;H^{l+1})\cap L^p_{\varepsilon}(I;W^{l,p}(\Omega))$, where $l\leq s$ and $k\leq r+1$. 
Then, the following convergence rates hold:
\begin{equation}\label{eq: L^2 space-time error}
\|(w-\mathcal{J}_{\tau,h}[w])_{tt}\|_{L^2_{\varepsilon}(L^2)}
\leq ch^{l+1}|w_{tt}|_{L^{2}_\varepsilon(H^{l+1})}+2c{C_{r,k}D_2\left(\frac{\tau}{\varepsilon}\right)}\tau^{k}|w|_{H^{k+2}_{\varepsilon}(L^2)},
\end{equation}
\begin{equation}\label{eq: H^1 space-time error}
\|\nabla (w-\mathcal{J}_{\tau,h}[w])\|_{L^2_{\varepsilon}(L^2)}
\leq ch^{l}| w|_{L^{2}_\varepsilon(H^{l+1})}+ 2c{C_{r,k}D_2\left(\frac{\tau}{\varepsilon}\right)}\frac{T^2}{\sqrt{2}}\tau^{k}|w|_{H^{k+2}_{\varepsilon}(H^1)},
\end{equation}
\begin{equation}\label{eq: L^p space-time error}
\|w-\mathcal{J}_{\tau,h}[w]\|_{L^p_{\varepsilon}(L^p)}
\leq ch^{l+1}|w|_{L^{p}_\varepsilon(W^{l,p})}+2c{C_{r,k}D_p\left(\frac{\tau}{\varepsilon}\right)}\frac{T^{2}}{2^{\frac{1}{p}}}\tau^{k}| w|_{W^{k+2,p}_{\varepsilon}(L^p)}.
\end{equation}
\end{lemma}
\begin{proof}

We decompose the approximation error as
\begin{align*}
    w - \mathcal{J}_{\tau,h}[w] = (w - \Pi_h w)+\Pi_h(w-\mathcal{J}_{\tau}[w]),
\end{align*}
and estimate the two errors separately. The term $w - \Pi_h w$ can be estimated directly using \eqref{eq: SZ error}. For $\Pi_h(w-\mathcal{J}_{\tau}[w])$, we use the stability of $\Pi_h$, \eqref{SZ_stab}, to deduce
\[
\|\Pi_h(w-\mathcal{J}_{\tau}[w])\|_{L^p_{\varepsilon}(X)}\le  2c \|w-\mathcal{J}_{\tau}[w]\|_{L^p_{\varepsilon}(X)}.
\]
Now, as $w,\mathcal{J}_{\tau}[w]\in \mathcal{K}(0,0)$, we apply \eqref{eq: second deriv norm Lk} with $q=p$ to conclude
\begin{equation}
\begin{aligned}
    \|w-\mathcal{J}_{\tau}[w]\|_{L^p_{\varepsilon}(L^p)} 
    \leq&
    \ 
    2^{-\frac{1}{p}}T^2 \|w_{tt} -\mathcal{I}_{\tau}^r[w_{tt}]\|_{L^p_{\varepsilon}(L^p)}
    \le   2^{-\frac{1}{p}}T^2{C_{r,k}D_{p}\left(\frac{\tau}{\varepsilon}\right)}\tau^{k} |w_{tt}|_{W^{k,p}_{\varepsilon}(L^p)},
    \label{eq: difference second deriv bound}
\end{aligned}
\end{equation}
from \eqref{eq: time Lp error}. Similarly employing \eqref{eq: second deriv norm}, one has
\begin{align*}
    \|\nabla (w-\mathcal{J}_{\tau}[w])\|_{L^2_{\varepsilon}(L^2)} \leq 2^{-\frac{1}{2}}T^2{C_{r,k}D_{2}\left(\frac{\varepsilon}{\tau}\right)}\tau^k|\nabla w_{tt}|_{H^{k}_{\varepsilon}(L^2)}.
\end{align*}
Also from \eqref{eq: time Lp error}, one has,
\begin{align*}
    \| w_{tt}- (\mathcal{J}_{\tau}[w])_{tt}\|_{L^2_{\varepsilon}(L^2)}  &=\| w_{tt}-\mathcal{I}_{\tau}^r[w_{tt}]\|_{L^2_{\varepsilon}(L^2)} \leq {C_{r,k}D_2\left(\frac{\tau}{\varepsilon}\right)}\tau^k|w_{tt}|_{H^k_{\varepsilon}(L^2)}.
\end{align*}
The above bounds along with triangle inequality and \eqref{eq: SZ error} are sufficient to show the required bounds.

\end{proof}

 We are now ready to prove the convergence of the finite element method.
\begin{theorem}[Convergence of the finite element method]\label{thm: convergence}
{With $p>2$}, let $u^{\varepsilon}$ and $u^{\varepsilon}_{\tau,h}$ be the solutions of \eqref{eq: Continuous Euler-lagrange equation} and \eqref{eq: Discrete Euler-lagrange equation} respectively and assume that $u^{\varepsilon}\in H^{k+2}_{\varepsilon}(I;H^1(\Omega))\cap W^{k+2,p}_{\varepsilon}(I;L^p(\Omega))\cap H^2_{\varepsilon}(I;H^{l+1})\cap L^p_{\varepsilon}(I;W^{l,p}(\Omega))$. Then, we have
\[
    \|\error\|_{E_{\varepsilon}}^2 +\alpha \|\error\|_{L^p_{\varepsilon}(L^p)}^p 
    \le C_4\Big( \tau^{\frac{kp}{p-1}} + h^{\min\{2l,\frac{(l+1)p}{p-1}\}}\Big),
\]
with $C_4>0$ depending on $p,T,\Omega,u^\varepsilon,f,s,r,D_{2}\left(\frac{\tau}{\varepsilon}\right),D_{p}\left(\frac{\tau}{\varepsilon}\right)$ and the shape-regularity of the mesh family.
\end{theorem}
\begin{proof}
Starting from Theorem \ref{theorem:ceas lemma}, setting $v=\mathcal{J}_{\tau,h}[u^{\varepsilon}]$, and using the best approximation estimates from Lemma \ref{lemma: projection error}, we deduce
\begin{align*}
    &\|\error\|_{E_{\varepsilon}}^2 +\alpha\|\error\|_{L^p_{\varepsilon}(L^p)}^p \\
    \leq& \  C\varepsilon^2 \Big(h^{2l+2}|u^{\varepsilon}_{tt}|_{L^{2}_\varepsilon(H^{l+1})}^2+\tau^{2k}|u^{\varepsilon}|_{H^{k+2}_{\varepsilon}(L^2)}^2\Big)
    +C\Big( h^{2l}|u^{\varepsilon}|_{L^{2}_\varepsilon(H^{l+1})}^2+\tau^{2k}| u^{\varepsilon}|_{H^{k+2}_{\varepsilon}(H^1)}^2\Big)\\ 
    &+\beta C\Big(h^{\frac{(l+1)p}{p-1}}|u^{\varepsilon}|_{L^{p}_\varepsilon(W^{l,p})}^{p'}+\tau^{\frac{kp}{p-1}}|u^{\varepsilon}|_{W^{k+2,p}_{\varepsilon}(L^p)}^{p'}\Big)
 +\gamma C\Big(h^{(l+1)p}|u^{\varepsilon}|_{L^{p}_\varepsilon(W^{l,p})}^p+\tau^{kp}| u^{\varepsilon}|_{W^{k+2,p}_{\varepsilon}(L^p)}^p\Big),
\end{align*}
from which the result follows.
\end{proof}

$C_4$ above depends on $\varepsilon$ in two places. Firstly the dependence on $D_{2}\left(\frac{\tau}{\varepsilon}\right),D_{p}\left(\frac{\tau}{\varepsilon}\right)$. For a fixed timestep $\tau$, these functions blow up exponentially as $\varepsilon \rightarrow 0$, with rate $e^{\frac{1}{\varepsilon}}$. However, for a given $\varepsilon$ there is a timestep $\tau_0$ such that $D_{2}\left(\frac{\tau}{\varepsilon}\right),D_{p}\left(\frac{\tau}{\varepsilon}\right)$ are bounded above uniformly independent of $\varepsilon$ for all $\tau \leq \tau_0$ sufficiently small. This in effect places a restriction on the timestep of the form $\tau \leq \mu\varepsilon$, for a constant $\mu>0$, such that these functions are bounded above.

The other dependence on $\varepsilon$ comes from the dependence of the higher derivatives of $u^{\varepsilon}$ in the weighted norms. For this exposition, we assume sufficient regularity of $u^{\varepsilon}$, say $u^{\varepsilon} \in C^{4}([0,T];C^2(\overline{{\Omega}}))$ such that the PDE  \eqref{eq: fourth order time PDE} with zero forcing and without nonlinearity but with nontrivial initial data holds pointwise. Using an asymptotic expansion to order $\mathcal{O}(\varepsilon^4)$, i.e. writing $u^{\varepsilon} = u_{comp}+ \mathcal{O}(\varepsilon^4)$, one can show that the higher derivatives of $u_{comp}$ in the $\|\cdot\|_{L^2}$ norm are not bounded above independent of $\varepsilon$, however they are bounded above independent of $\varepsilon$ in the weighted $\|\cdot\|_{L^2_{\varepsilon}}$ norms, which indicates that $u^{\varepsilon}$ will also be bounded independent of $\varepsilon$ in these weighted norms. The full details of the asymptotic expansions are given in Appendix \ref{sec: Asymptotic analysis}.

For weak nonlinearities we also have the following result, whose proof is completely analogous and, thus, omitted for brevity. 
\begin{corollary}[Optimal rates for small nonlinearities]\label{cor: small p convergence}
If the assumptions of Theorem \ref{cor: small nonlin cea} hold, then
\begin{equation*}
    \|\error\|_{E_{\varepsilon}} \le  \min\left\{(1+G)C_5\big( \tau^k + h^{l}\big),\sqrt{C_4}\left(\tau^{\frac{kp}{2(p-1)}}+h^{\min\{l,\frac{(l+1)p}{2(p-1)}\}}\right)\right\},
\end{equation*}
with $G =\mathcal{O}(e^{\frac{T}{\varepsilon}})$ and $C_5$ depending on $p,T,\Omega,u^{\varepsilon},f,s,r,D_2\left(\frac{\tau}{\varepsilon}\right)$ and the shape-regularity of the mesh family.\qed
\end{corollary}
\begin{remark}
For the linear Klein-Gordon equation, corresponding to $p=2$, the above imply the convergence rate $\|e_{\varepsilon}\|_{KG}\leq C_6(\tau^k+h^l)$, with $C_6$ having the same dependencies as $C_5$. For the linear wave equation, corresponding to $p=0$, the respective convergence rate reads $\|e_{\varepsilon}\|_{E_{\varepsilon}}\leq C_{5}(\tau^k+h^l)$. In particular, there is no dependence on $e^{\frac{T}{\varepsilon}}$ in either $C_5$ or $C_6$.
\end{remark}

\begin{remark} The results of Theorem \ref{thm: convergence} and Corollary \ref{cor: small p convergence} also apply to the semilinear Klein-Gordon, i.e. a reaction term in \eqref{eq: Model equation} of the form $u+\frac{p}{2}|u|^{p-2}u$ under the same conditions, by simply substituting the energy norm $\|\cdot\|_{E_{\varepsilon}}$ with $\|\cdot\|_{KG}$. This follows by the inclusion of an additional linear reaction term $u$, which affects only the linear part of the spatial operator.
    
\end{remark}

\section{Numerical experiments}
\label{sec: numerical experiments}
\subsection{Temporal discretisation}
For the time discretisation, we use the space of maximal regularity cubic splines $\mathcal{S}^{\Delta_{\tau}}(3)$, whereby $\text{dim}(\mathcal{S}^{\Delta_{\tau}}(3)) = N_t+3$. The basis is defined via the classical B-spline construction (see, e.g., \cite[Chapter~11]{MR2006500}), which have compact support. For simplicity, we consider uniform timestepping, allowing for an explicit representation of the basis functions: $\mathcal{S}^{\Delta_{\tau}}(3) =\text{span}\{\varphi_{-1},\varphi_{0},\cdots,\varphi_{N_t},\varphi_{N_t+1}\}$, where
\begin{equation}
    \varphi_k(t) = \frac{1}{4\tau^3}\Phi(t-(k-2)\tau),
    \label{eq: basis function definition}
\end{equation}
and
$    \Phi(t+2\tau) = (t+2\tau)_{+}^3-4(t+\tau)_{+}^3+6(t)_{+}^3-4(t-\tau)_{+}^3+(t-2\tau)_{+}^3$,
with $(x)_{+} = x \mathbf{1}_{x\geq 0}$. We have $\text{supp}(\varphi_k)=[(k-2)\tau,(k+2)\tau]$ and the normalisation $\varphi_k(k\tau)=1$. We stress that non-uniform timestepping is by all means possible using the classical De Boor's algorithm. 

To enforce zero initial position and velocity, we modify the basis of $\mathcal{S}^{\Delta_{\tau}}(3)$ to contain the functions ${1,t}$ instead: $\varphi_2,\cdots,\varphi_{N_t+1}$ are all linearly independent of ${1,t}$ owing to their compact support not intersecting the origin. Since $\mathrm{dim}\{1,t,\varphi_2,\cdots,\varphi_{N_t+1}\} = N_t+2$, we require an additional basis function obtained as a linear combination of the three B-splines that take nonzero values at $t=0$: writing $\tilde{\varphi} = a\varphi_{-1}+b\varphi_0+c\varphi_1$, we enforce $\tilde{\varphi}(0)=0$ and $ \tilde{\varphi}_t(0)=0$, giving $a\varphi_{-1}(0)+b\varphi_0(0)+c\varphi_1(0)=0$ and $a(\varphi_{-1})_t(0)+b(\varphi_{0})_t(0)+c(\varphi_{1})_t(0) = 0$, respectively.
We also impose the arbitrary condition $\tilde{\varphi}(\tau)=1$. Solving these three equations gives
\begin{align*}
    \tilde{\varphi} = \frac{8}{7}\varphi_{-1}-\frac{4}{7}\varphi_0+\frac{8}{7}\varphi_1.
\end{align*}
Hence, $\mathcal{S}^{\Delta_{\tau}}(3)=\text{span}\{1,t,\tilde{\varphi},\varphi_2,\cdots,\varphi_{N_t},\varphi_{N_t+1}\}$. We solve on $U_{\tau} =\mathrm{span}\{\tilde{\varphi},\varphi_2,\cdots,\varphi_{N_t},\varphi_{N_t+1}\}$, satisfying the initial conditions.
For convenience, we now relabel these functions as $\phi_1 = \tilde{\varphi},\phi_i = \varphi_i$, $i=2,\cdots,N_t+1$ so that
\[
U_{\tau} =\mathrm{span}\{\phi_1,\phi_2,\cdots,\phi_{N_t+1}\},
\]
keeping in mind that $\phi_1$ is not a translation of any of the $\phi_2,\cdots,\phi_{N_t+1}$.

\subsection{Linear problem}
\subsubsection{Zero-dimensional case}
To illustrate the challenges, we begin by studying the ODE problem
\begin{equation}
\begin{split}
    &\varepsilon^2u_{tttt}^{\varepsilon}-2\varepsilon u_{ttt}^{\varepsilon}+u^{\varepsilon}_{tt}+\lambda u^{\varepsilon} = f, \\
    & u^{\varepsilon}(0)=0,u_t^{\varepsilon}(0)=0,\quad 
    u^{\varepsilon}_{tt}(T)=u^{\varepsilon}_{ttt}(T)=0.
\end{split}
\label{eq: fourth order ODE}
\end{equation}
This problem is relevant as it can be obtained by expanding the solution of \eqref{eq: fourth order time PDE} in a series of eigenfunctions of the Laplacian, with corresponding eigenvalues $\lambda> 0$. The method corresponding to \eqref{eq: fourth order ODE} reads: find $u^{\varepsilon}_{\tau} \in U_{\tau}$, such that
\begin{equation}
    \int_0^Te^{-\frac{t}{\varepsilon}}\left(\varepsilon^2(u^{\varepsilon}_{\tau})_{tt}v_{tt}+\lambda u^{\varepsilon}_{\tau}v\right)\;\text{d}t = \int_0^Te^{-\frac{t}{\varepsilon}}fv\;\text{d}t, \quad \text{for all } v \in U_{\tau}.
\end{equation}
This reduces to solving the linear system
\begin{equation}
    \left(\varepsilon^2K_{\tau}+\lambda L_{\tau}\right)\sigma = F,
    \label{eq: 0D matrix equation}
\end{equation}
where $\sigma=(\sigma_j)_{j=1}^{N_t+1}$, $u^{\varepsilon}_{\tau}= \sum_{j=1}^{N_t+1}\sigma_j\phi_j(t)$, $L_{\tau}=(l_{i,j})_{i,j=1}^{N_t+1}$ and $K_{\tau}=(k_{i,j})_{i,j=1}^{N_t+1}$ are the mass and stiffness matrices, $ F = (f_{i})_{i=1}^{N_t+1}$ the load vector, where
\begin{align*}
    l_{i,j}=\int_0^Te^{-\frac{t}{\varepsilon}}\phi_j\phi_i\,\ud t,\;\; k_{i,j}=\int_0^Te^{-\frac{t}{\varepsilon}}\phi_{j,tt}\phi_{i,tt}\,\ud t, \;\; f_{i} =\int_0^Te^{-\frac{t}{\varepsilon}}f(t)\phi_i(t)\,\ud t.
\end{align*}
Owing to the ellipticity of the underlying problem,  the system matrix $\varepsilon^2 K_\tau + \lambda L_\tau$ is symmetric and positive definite, and hence a unique solution to \eqref{eq: 0D matrix equation} exists. However, the matrix entries need to be computed carefully and the system is badly conditioned as we see next.

\subsubsection{Underflow errors}\label{subsec: underflow}
A naive implementation of \eqref{eq: 0D matrix equation} will quickly fail due to the weighted inner products in the definition of $L_{\tau},K_{\tau}$. For example, considering $\varepsilon = 0.002$, $T=2$, computation of  $e^{-2000}\approx 10^{-1000}$ which leads to an underflow error. Considering these integrals individually, for $i,j>1$ we can exploit the compact support of B-splines to write
\begin{equation}
\begin{split}
l_{i,j} &= \int_0^Te^{-\frac{t}{\varepsilon}}\phi_j(t)\phi_i(t)\;\text{d}t \stackrel{\text{B-spline}}{=} \int_{\max\{0,(i-2)\tau\}}^{\min\{T,(i+2)\tau\}}e^{-\frac{t}{\varepsilon}}\phi_j(t)\phi_i(t)\;\text{d}t \\
    &= \frac{1}{16\tau^6}\int_{\max\{0,(i-2)\tau\}}^{\min\{T,(i+2)\tau\}}e^{-\frac{t}{\varepsilon}}\Phi(t-(j-2)\tau)\Phi(t-(i-2)\tau)\;\text{d}t \\
    &=\frac{e^{-\frac{i\tau}{\varepsilon}}}{16\tau^6}\int_{\max\{-i\tau,-2\tau\}}^{\min\{T-i\tau,2\tau\}}e^{-\frac{s}{\varepsilon}}\Phi(s-((j-i)-2)\tau)\Phi(s-(-2\tau))\;\text{d}s \\
    &= e^{-\frac{i\tau}{\varepsilon}}\int_{\max\{-i\tau,-2\tau\}}^{\min\{T-i\tau,2\tau\}}e^{-\frac{s}{\varepsilon}}\phi_{j-i}(s)\phi_0(s)\;\text{d}s,
\end{split}
\label{eq: exp shift}
\end{equation}
where  $\phi_0 = \varphi_0$ and $\phi_{j-i}=\varphi_{j-i}$ for $j-i\leq 0$. A similar calculation can be performed for $F$ and $K_{\tau}$. Defining
\[
    \tilde{l}_{i,j} = e^{\frac{i\tau}{\varepsilon}}l_{i,j}, \quad
    \tilde{k}_{i,j} = e^{\frac{i\tau}{\varepsilon}}k_{i,j},\quad
    \tilde{f}_{i} = e^{\frac{i\tau}{\varepsilon}}f_{i},
\]
we see that $\sigma$ is also the solution of the system
\begin{equation}
    \left(\varepsilon^2 \tilde{K}_{\tau}+\lambda \tilde{L}_{\tau}\right)\sigma = \tilde{F}.
\label{eq: exp linear system}
\end{equation}
Let $E = \mathrm{diag}(e^{\frac{\tau}{\varepsilon}},e^{\frac{2\tau}{\varepsilon}},\cdots,e^{\frac{T}{\varepsilon}},e^{\frac{T+\tau}{\varepsilon}})$, then $EK_{\tau} = \tilde{K}_{\tau},EL_{\tau} = \tilde{L}_{\tau}$. In particular, in the case $i=1$ or $j=1$, which corresponds to the basis function $\phi_1 = \tilde{\varphi}$, $\tilde{L}_{\tau},\tilde{K}_{\tau}$ is defined so that this relation holds. 

In practice, we solve \eqref{eq: exp linear system} instead of \eqref{eq: 0D matrix equation}. A quadrature rule could be devised for computing these integrals exactly for polynomials, however we find that standard quadrature rules work well in practice \cite{MR712135,2020SciPy-NMeth}.

\begin{remark}
The above process is equivalent to preconditioning the system by the matrix $E$ and is also equivalent to considering $\{\phi_i\}_{i=1}^{N_t+1}$ as the basis for the trial space and $\{e^{\frac{i\tau}{\varepsilon}}\phi_i\}_{i=1}^{N_t+1}$ as the basis for the test space.
\end{remark}

\subsubsection{Ill-conditioning}
\begin{figure}
    \centering
    \includegraphics[width=0.6\linewidth]{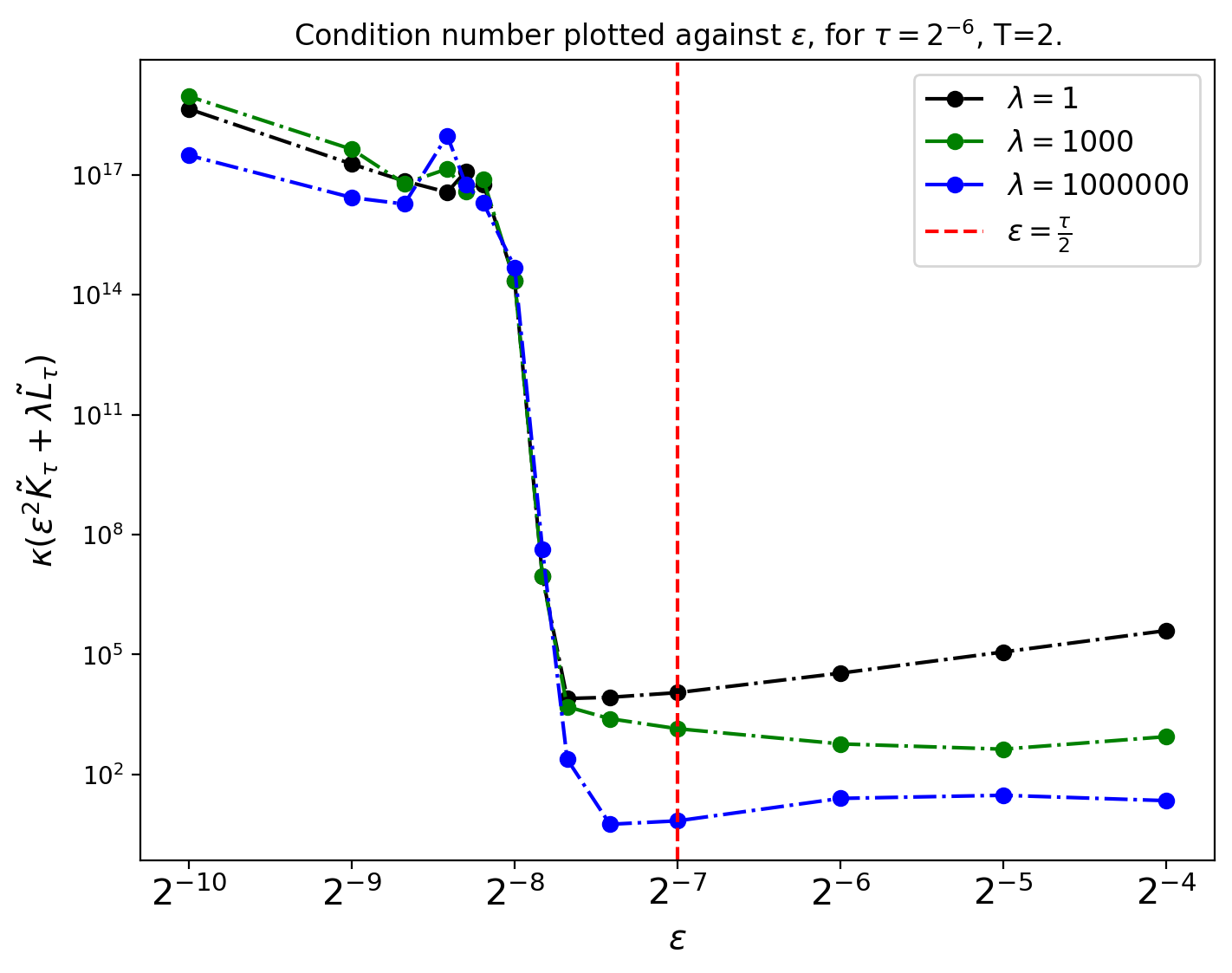}
    \caption{ The condition number of the matrix in \eqref{eq: exp linear system}, depicting $\kappa(\varepsilon^2\tilde{K}_{\tau}+\lambda \tilde{M}_{\tau})$ against $\varepsilon$ for fixed $\tau=2^{-6}$ with $T=2$, for various values of $\lambda =1,1000,1000000$.}
    \label{fig:Condition number plot}
\end{figure}

As perhaps expected, the matrix in \eqref{eq: exp linear system} becomes increasingly ill-conditioned as $\varepsilon\rightarrow 0^+$. In practice, we see that the condition number grows when $\varepsilon< \tau$. The dependence of the condition number on the parameter $\varepsilon$ for fixed $\tau$ is illustrated in Figure \ref{fig:Condition number plot}. In particular, we note that conditioning depends chiefly on the parameters $\varepsilon$ and $\tau$, and only mildly on $\lambda$. The problem appears to be well conditioned for $\varepsilon>\frac{\tau}{2}$. 

It appears that ill-conditioning is due to the weight $e^{-\frac{s}{\varepsilon}}$. Indeed, for an entry \[
\tilde{l}_{i,j} = \int_{\max\{-i\tau,-2\tau\}}^{\min\{T-i\tau,2\tau\}}e^{-\frac{s}{\varepsilon}}\phi_{j-i}(s)\phi_0(s)\,\ud s\] of the mass matrix $\tilde{L}_{\tau}$, the weight $e^{-\frac{s}{\varepsilon}}$ takes values between $e^{-\frac{2\tau}{\varepsilon}}$ and $e^{\frac{2\tau}{\varepsilon}}$. For large ratios of $\frac{\tau}{\varepsilon}$, this will lead to the integrand being localised at $s=-2\tau$. We find that restricting the ratio of $\frac{\tau}{\varepsilon}$ to $\frac{\tau}{\varepsilon}\leq 2$ suffices to ensure that the problem is well conditioned, as this results in the weight taking values in $[e^{-4},e^{4}]$. 

\subsection{One-dimensional case}

Throughout, the domain is restricted to $\Omega=[0,1]$, which we split into equal subintervals of width $h$, such that $N_x = \frac{1}{h}$. The spatial discrete space consists of continuous, piecewise degree-one polynomials, i.e., $s=1$.

For the linear problem with trivial initial data, i.e.,  $u^{\varepsilon}_{\tau,h} \in \mathcal{K}_{\tau,h} = U_{\tau}\otimes V_h$, we consider the spatial mass and stiffness matrices $M_h$ and $A_h$, respectively. Then, the solution of the space-time method is given by solving the linear system
\begin{equation}
\begin{split}
 \left(\varepsilon^2K_{\tau}\otimes M_h+L_{\tau}\otimes A_h\right)\mathbf{\sigma} = F,
\end{split}    
\label{eq: Tensor product system}
\end{equation}
where $\otimes$ denotes the Kronecker product. In the calculation of the temporal inner products, the underflow-avoiding considerations from Section \ref{subsec: underflow} are used.
We begin by investigating the error $\error=u^{\varepsilon} -u^{\varepsilon}_{\tau,h}$. With the choice of the discretisation space $\mathcal{K}_{\tau,h}$, we expect second-order convergence with respect to $\tau$ and first order in $h$, as $s=1$. Measuring the error in the energy norm $E_{\varepsilon}$ hinders insightful information: the weight $e^{-\frac{t}{\varepsilon}}$ supresses any nontrivial errors away from $t=0$. Instead, we estimate the error in the $L^2(L^2)$ norm and the $H^1(L^2)$, $H^2(L^2)$ seminorms combined with the $L^2(H^1)$  seminorm.

\begin{figure}
    \centering
    \begin{subfigure}[b]{0.49\textwidth}
    \includegraphics[scale=0.44]{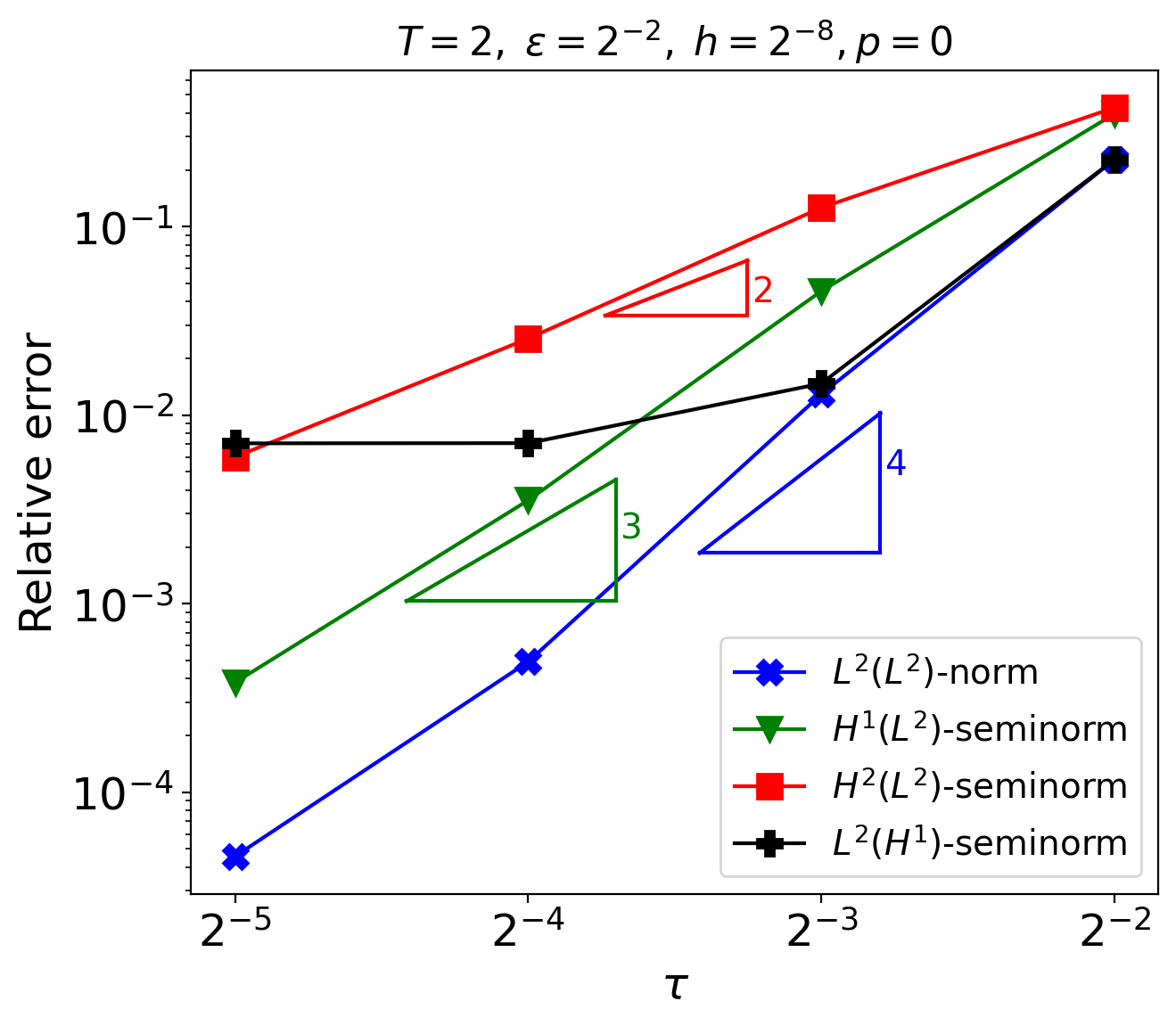}  
        \caption{}
        \label{fig: linear tau}
    \end{subfigure}
    \begin{subfigure}[b]{0.49\textwidth}
    \includegraphics[scale=0.44]{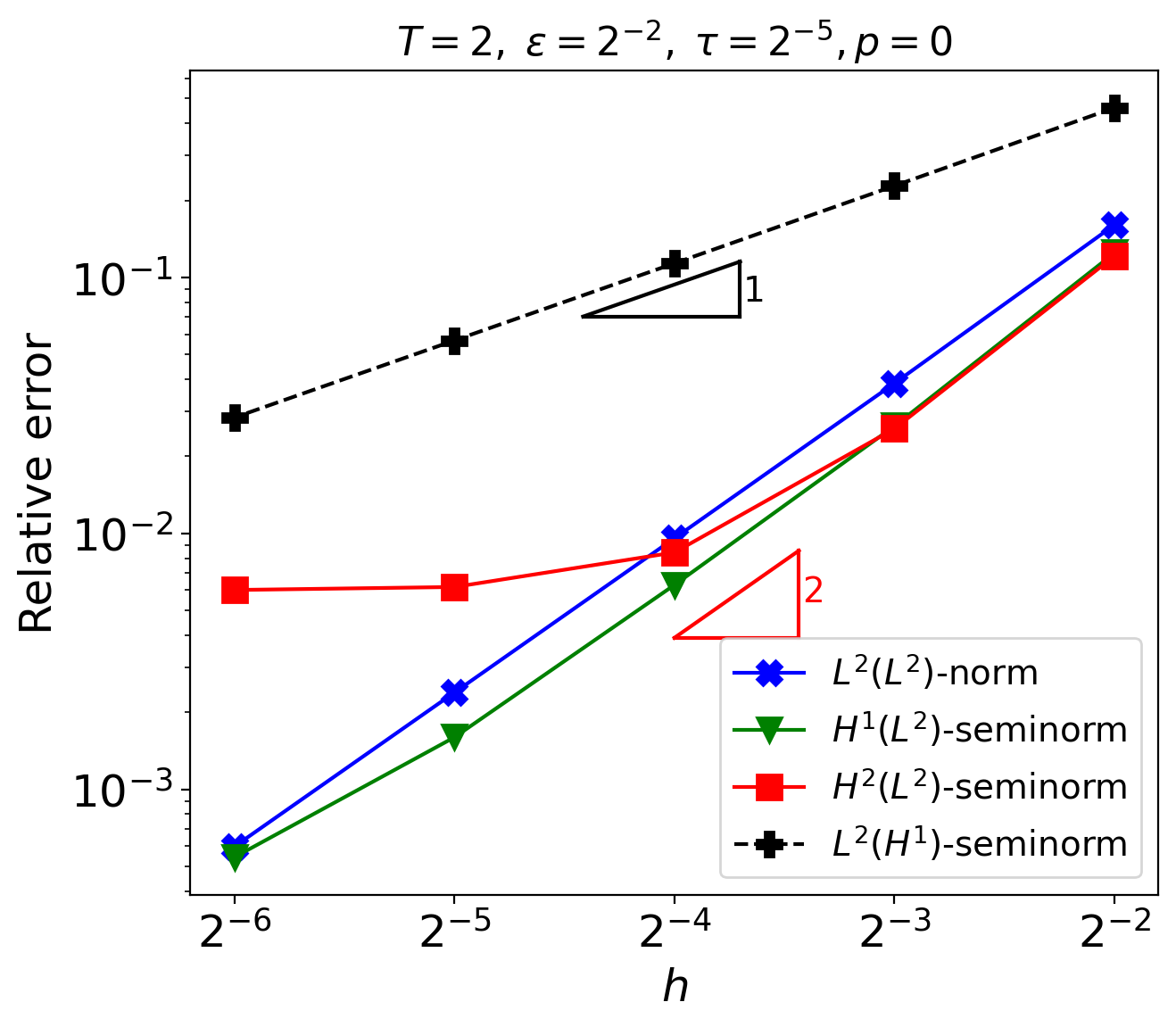}  
        \caption{}
        \label{fig: linear space}
    \end{subfigure}
\caption{Convergence history in the parameters $\tau,h$ to the linear, regularised  problem with solution $u^{\varepsilon}$, for $\varepsilon = 2^{-2}$ fixed. In \ref{fig: linear tau}, $h=2^{-8}$ is fixed with $\tau$ decreasing. In \ref{fig: linear space}, $\tau=2^{-5}$,is fixed with decreasing $h$.} 
\label{fig: linear convergence rates}
\end{figure}

In the test case, we will set $T=2$, $\Omega=[0,1]$, $\varepsilon = 0.25$ and choose the forcing $f_\varepsilon(x,t)$ such that $u^{\varepsilon}(x,t) = \sin(2\pi x)(T-t)^2\sin^2(2\pi t)$ solves the fourth order PDE model problem \eqref{eq: fourth order time PDE} for $d=1$, without any linear or nonlinear reaction term.
The computer convergence history is presented in Figure \ref{fig: linear convergence rates}.

\begin{remark}
The increased rate in Figure \ref{fig: linear convergence rates} for weaker-than-energy norms can be shown via an Aubin-Nitsche argument, upon assuming an elliptic regularity estimate for the dual problem of\eqref{eq: Continuous Euler-lagrange equation}.
\end{remark}

\subsection{Nonlinear problem}
For the nonlinear problem, we employ Newton iterations. 
Denote the approximation of $u^{\varepsilon}_{\tau,h}$ at the $k^{\text{th}}$ Newton iteration by $u^{\varepsilon,(k)}_{\tau,h}$. The update step is given by 
\begin{align*}
    u^{\varepsilon,(k+1)}_{\tau,h} = u^{\varepsilon,(k)}_{\tau,h}+\delta u^{(k)},
\end{align*}
where $\delta u^{(k)}$ satisfies
\begin{equation}
    a_u(u^{\varepsilon,(k)}_{\tau,h};v,\delta u^{(k)}) = -\left(A(u^{\varepsilon,(k)}_{\tau,h},v)+B(u^{\varepsilon,(k)}_{\tau,h};v)-f(v)\right), 
\end{equation}
for all $v \in \mathcal{K}_{\tau,h}$, with \begin{align*}
    &a_{u}(u;v,w) := A(v,w)+
    \frac{p(p-1)}{2}\int_0^Te^{-\frac{t}{\varepsilon}}\left<|u|^{p-2}w,v\right>_{L^2(\Omega)}\,\ud t.
\end{align*}
This formulation is derived by computing the Fréchet derivative of the weak formulation with respect to the first entry.
To set up the linear system, at each iteration integrals involving the nonlinear terms need to be computed. As the initial guess, $u_{\tau,h}^{\varepsilon,(0)}$, we use the solution of the linear problem \eqref{eq: Tensor product system}. In the calculation of the temporal inner products, the underflow-avoiding considerations from Section \ref{subsec: underflow} are again used. The stopping criterion $\frac{||\mathbf{\sigma}^{(k)}||_{l^2}}{\#\mathrm{dof}}\leq 10^{-10}$ is used, where $\mathbf{\sigma}^{(k)}$ is the vector of coefficients of $\delta u^{(k)}$ expanded in the discrete space-time basis.  

\begin{figure}
    \centering
    \begin{subfigure}[b]{0.49\textwidth}
    \includegraphics[scale=0.44]{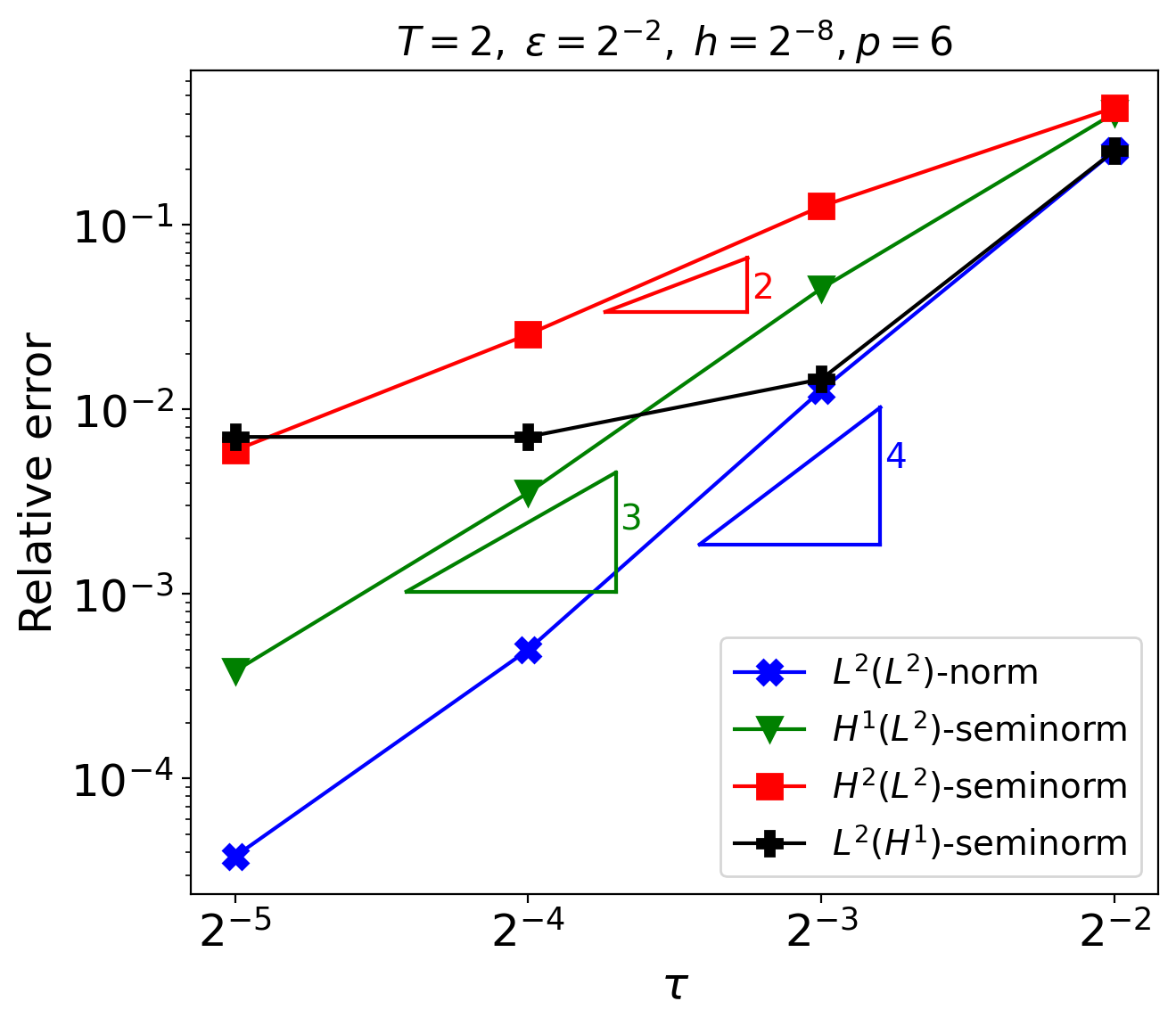}  
        \caption{}
        \label{fig: nonlinear tau}
    \end{subfigure}
    \begin{subfigure}[b]{0.49\textwidth}
    \includegraphics[scale=0.44]{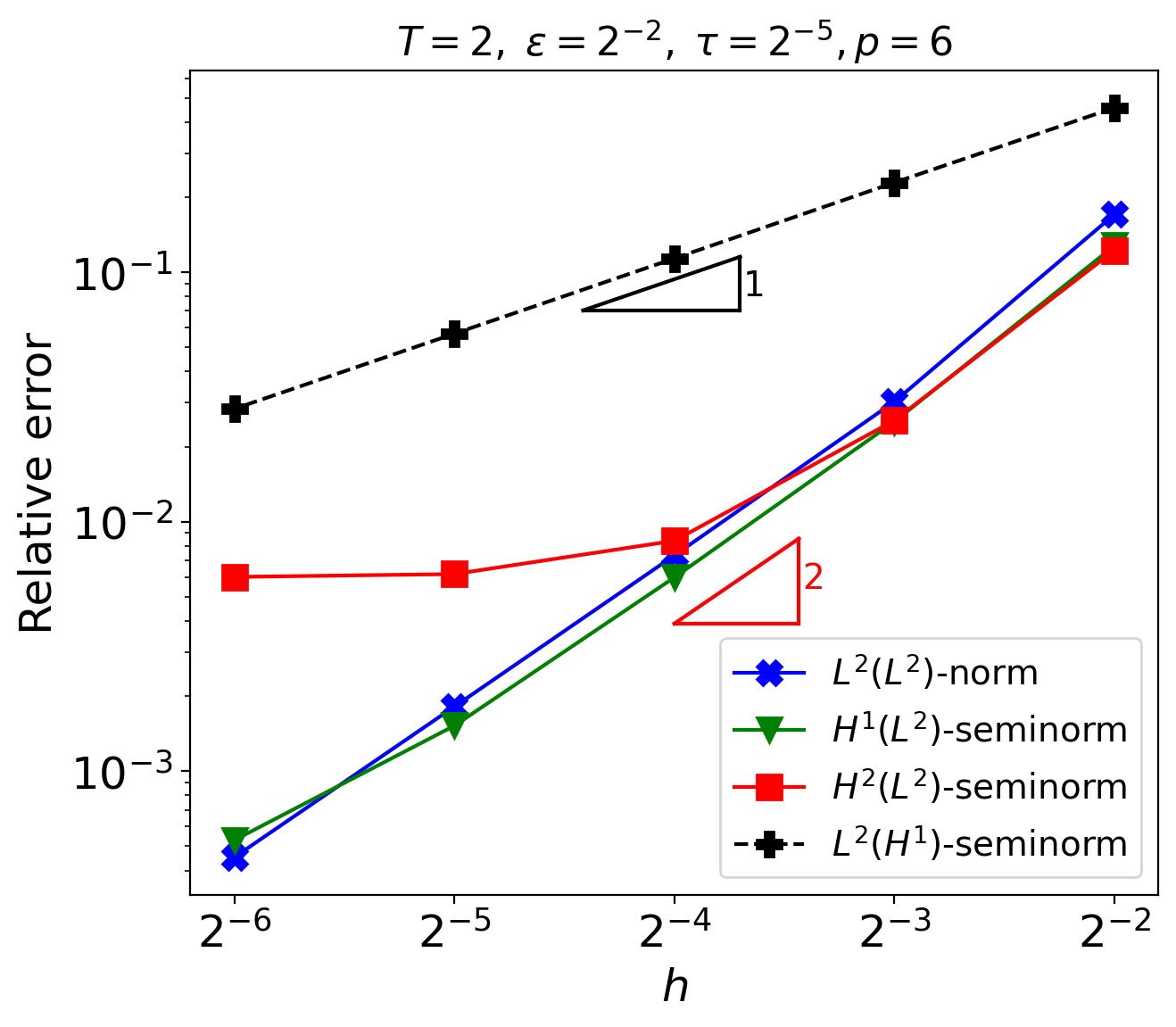}  
        \caption{}
        \label{fig: nonlinear space}
    \end{subfigure}
\caption{Two plots indicating the convergence rates in the parameters $\tau,h$ to the nonlinear regularised $u^{\varepsilon}$ problem with $\varepsilon = 2^{-2}$ fixed. In \ref{fig: nonlinear tau}, $h=2^{-8}$ is fixed with $\tau$ decreasing. In \ref{fig: nonlinear space}, $\tau=2^{-5}$,is fixed with decreasing $h$.}
\label{fig: nonlinear convergence rates}
\end{figure}

We set $p=6$ and, again, choose the forcing $f_{2,\varepsilon}$ such that $u^{\varepsilon}(x,t) = \sin(2\pi x) (T-t)^2 \sin^2(2\pi t)$.
The convergence rates are presented in Figure \ref{fig: nonlinear convergence rates}; they are as predicted by  Corollary~\ref{cor: small p convergence} for $d = 1$.

So far, we have discussed the convergence of the space-time method to the exact solution $u^\varepsilon$ of the \emph{perturbed} fourth-order problem. We now discuss the utility of the proposed method to approximate solution to the original wave problem \eqref{eq: Model equation}, that is, we assess the error $u^0-u_{\tau,h}^{\varepsilon}$, as $\tau,h,\varepsilon\to 0$. 

To ensure tractable condition numbers, we select $\tau \leq 2\varepsilon$ as suggested by Figure \ref{fig:Condition number plot} and, in particular, set $\varepsilon = \frac{\tau}{2}$. {This timestep restriction also ensures that $D_p\left(\frac{\tau}{\varepsilon}\right)\leq 2$ in Theorem \ref{thm: convergence}.} Since $|u^0-u^\varepsilon|=\mathcal{O}(\varepsilon)$, cf., Section \ref{BL}, this effectively places a total first order convergence rate. With $\tau, \varepsilon$ chosen this way, we may choose $h  = \sqrt{\tau}$ so that in the $L^2$-norm the errors are balanced. 

For $p=4$, we choose $f(x,t) = 2\pi \sin(\pi x)\cos (\pi t)+ 2\sin^3(\pi x)t^3\sin^3(\pi t)$ such that $u(x,t) = \sin(\pi x) t \sin(\pi t)$, $x\in(0,1)$, $t\in(0,2)$, is the exact solution to \eqref{eq: Model equation} with trivial initial data.

\begin{figure}
    \centering
    \includegraphics[width=0.55\linewidth]{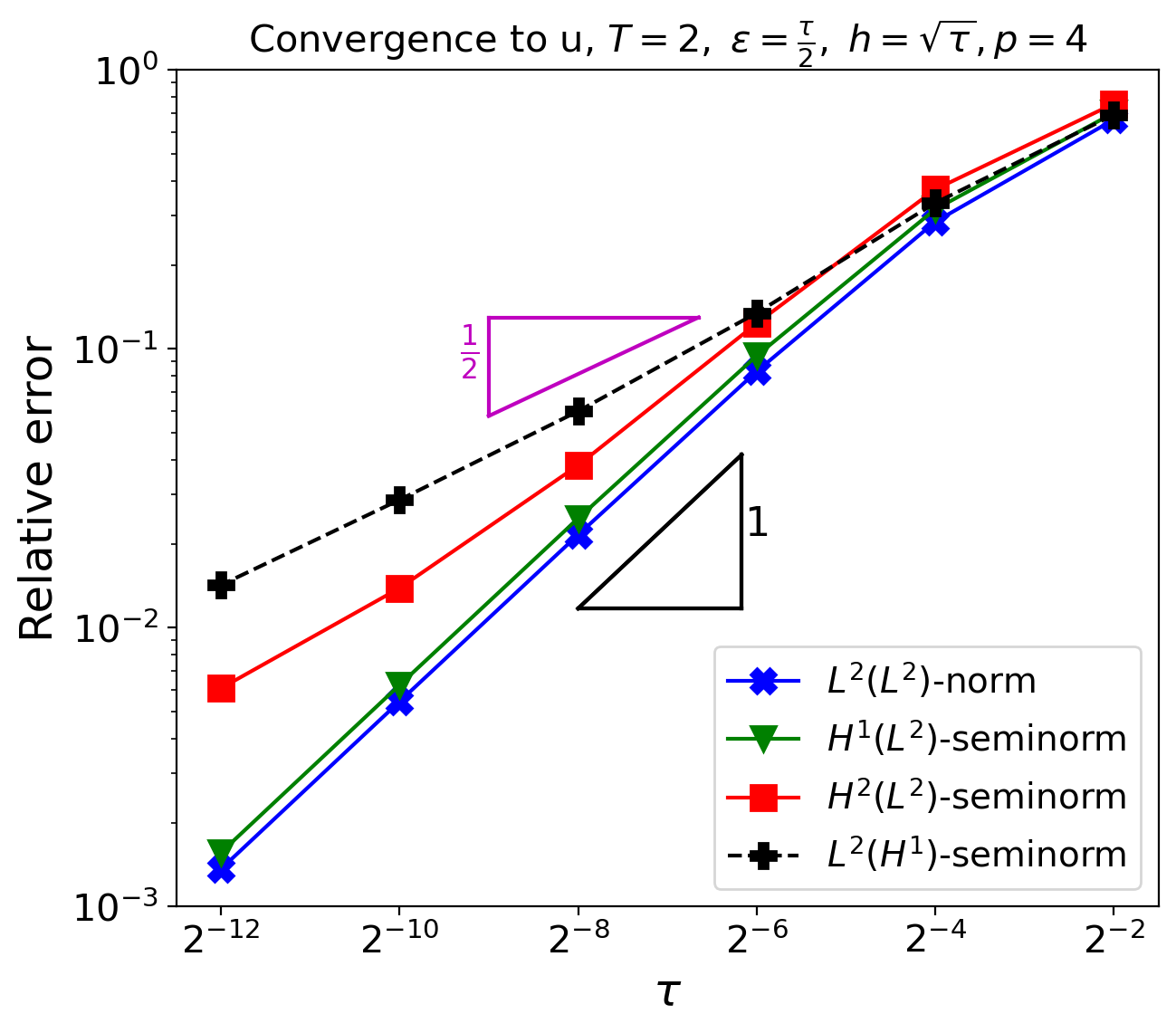}
    \caption{Convergence history for the nonlinear wave problem \eqref{eq: Model equation} problem with $p=4$. Here $h$ and $\varepsilon$ are chosen such that $h=\sqrt{\tau}$ and $\varepsilon = \frac{\tau}{2}$, with $T=2$.}
    \label{fig:convergence all together}
\end{figure}

The convergence rates of $u^{\varepsilon}_{\tau,h}\to u$ for $\varepsilon = \frac{\tau}{2}$ and $h =\sqrt{\tau}$ are presented in Figure \ref{fig:convergence all together}.

We note that, in an higher machine precision setting, or with more elaborate preconditioning strategies, an arbitrary order convergence scheme is possible, in principle, by selecting $\varepsilon = \tau^{r}$ for $r$ sufficiently large.

\section{Conclusions}\label{sec:conclusions}

Motivated by De Giorgi's elliptic reguralisation principle, this work investigates the construction of conforming, Galerkin-type numerical methods for the semilinear wave equation in second-order formulation. The underlying ellipticity gives well-posedness on the discrete level for any conforming discretisation, thereby guaranteeing  unconditional stability of the method. 

The Galerkin projection, along with the monotonicity of the nonlinearity, allow for the proof of quasi-optimality. \emph{A priori} bounds are also shown for sufficiently smooth solutions. In the cases of the energy norm controlling the nonlinearity, we are able to prove that there is no loss in the convergence rate.

The challenges of numerically solving the problem are illustrated through the zero-dimensional ODE setting. In particular, the underflow error due to the exponential weight is addressed by a basic preconditioning strategy. A study of the connection between the ratio of the regularisation parameter $\varepsilon$ and the timestep $\tau$ to the condition number for the resulting systems is performed to assess the practicality of the approach as a viable numerical method. A series of numerical experiments in a one-dimensional setting confirm the theoretical findings.

The motivation for the present work has been the use of elliptic reguralisation in the context of rigorous numerical methods for nonlinear wave equations with view to tackling in the future highly nonlinear complex wave problems \cite{MR4913691,MR4156784,WIDE} that currently lack satisfactory numerical treatment. This work is a first step in this challenging direction. 

\section{Acknowledgements}
\noindent BH was supported by the EPSRC Centre for Doctoral Training in Mathematical Modelling, Analysis and Computation (MAC-MIGS), funded by the UK Engineering and Physical Sciences Research Council (grant EP/S023291/1, Heriot-Watt University and the University of Edinburgh). EHG gratefully acknowledges the financial support of EPSRC (grant number EP/W005840/2).
\bibliographystyle{amsplain}
\bibliography{bibliography}
\appendix
\section{Asymptotic analysis of the boundary layer}
\label{sec: Asymptotic analysis}
We study the asymptotic behaviour of the boundary value problem in time, highlighting the relevance of the exponentially weighted norms. For simplicity, we restrict to the linear boundary value problem ($p=0$) for $\varepsilon\ll 1$:
\begin{equation}
\begin{split}
\varepsilon^2u^{\varepsilon}_{tttt}-2\varepsilon u^{\varepsilon}_{ttt}+u^{\varepsilon}_{tt}-\Delta u^{\varepsilon}=0&\text{ in } \Omega \times (0,T), \\
u^{\varepsilon}=0, &\text{ on }\partial\Omega \times[0,T], \\
u^{\varepsilon}(\cdot,0)=u^0(\cdot),u_t^{\varepsilon}(\cdot,0)=u^1(\cdot), &\text{ in }\Omega \\
\varepsilon u^{\varepsilon}_{tt}(\cdot,T)=0,\varepsilon^2 u_{ttt}^{\varepsilon}(\cdot,T)=0&\text{ in } \Omega.
\end{split}
\label{eq: fourth order time PDE given IC}
\end{equation}
\subsection{Outer region}
We begin with the outer problem and seek a solution to \eqref{eq: fourth order time PDE given IC} to order $\mathcal{O}(\varepsilon^4)$. We write
\begin{align*}
    u^{\varepsilon}(x,t) = u_0(x,t)+\varepsilon u_1(x,t)+\varepsilon^2u_2(x,t)+\varepsilon^3u_3(x,t)+\mathcal{O}(\varepsilon^4).
\end{align*}
Then, $u_0,u_1,u_2,u_3$ satisfy the problems:\\
\begin{minipage}{0.49\textwidth}
\begin{align*}
    &u_{0,tt}-\Delta u_0 = 0  \text{ in }\Omega \times I, \\
    &u_0 =  0 \text{ on }\partial\Omega \times I, \\
    &u_0(\cdot,0) = u^0(\cdot),\  u_{0,t}(\cdot,0)= u^1(\cdot) \ \text{ in }\Omega, 
\end{align*}
\end{minipage}
\begin{minipage}{0.45\textwidth}
\begin{align*}
   & u_{m,tt}-\Delta u_m = 2 u_{m-1,ttt}+u_{m-2,tttt}& \text{ in }\Omega \times I, \\
    &u_m =0 \text{ on }\partial\Omega \times I, \\
    &u_m(\cdot,0) = 0=u_{m,t}(\cdot,0) \text{ in }\Omega,
\end{align*}
\end{minipage}\smallskip\\
respectively, for $m=1,2,3$ with the convention that $u_{-1}=u_{-2}=0$.
Thus, $u_m$ can be computed in a sequential manner, and we denote the outer solution as:
   $ u_{out} = u_0+\varepsilon u_1+\varepsilon^2u_2+\varepsilon^3u_3$.
\subsection{Inner region}
For the inner region, we denote $U(s) = u^{\varepsilon}(T-\varepsilon s)$. Let $'$ denote differentiation with respect to the variable $s$. Substituting into \eqref{eq: fourth order time PDE given IC}, we obtain the problem for the inner region: 
\begin{equation}
\begin{split}
    U''''+2U'''+U'' +\varepsilon^2\Delta U = 0,\qquad   U''(0)=0,\ U'''(0)=0.
\end{split}
\label{eq: inner region pde}
\end{equation}
We seek a solution to \eqref{eq: inner region pde} again to order $\mathcal{O}(\varepsilon^4)$, and we write
 $   U = U_0+\varepsilon U_1+\varepsilon^2 U_2+\varepsilon^3 U_3+\mathcal{O}(\varepsilon^4)$.
Then, for $m=1,2,3,4$, $U_m$ satisfies
\begin{align*}
    U_{m}''''+2U_{m}'''+U_m''=-\Delta U_{m-2}, \qquad U_{m}''(0)=0,\ U_{m}'''(0)=0,
\end{align*}
with the convention $U_{-1}=U_{-2}=0$. Hence,
\begin{align*}
    U_0(x,s) = A_0(x)+B_0(x) s, \qquad  
    U_1(x,s) = A_1(x)+B_1(x)s,
\end{align*}
with $A_{0},A_1,B_0,B_1$ to be determined later by matching. $U_2$ and $U_3$ satisfy, in turn,
\begin{align*}
    U_{2}''''+2U_2'''+U_2'' = -(\Delta A_0+\Delta B_0 s),\qquad   U_2''(0)=U_2'''(0)=0,\\
    U_{3}''''+2U_3'''+U_3'' = -(\Delta A_1+\Delta B_1 s),\qquad 
    U_3''(0)=U_3'''(0)=0,
\end{align*}
whose exact solutions are given by
\begin{align*}
    U_i(x,s) =&\  A_i+B_i s-\frac{\Delta B_{i-2}}{6}s^3+\Big(\Delta B_{i-2}-\frac{\Delta A_{i-2}}{2}\Big)s^2\\
    &+\left[\left(-4\Delta B_{i-2}+3\Delta A_{i-2}\right)+(\Delta A_{i-2}-\Delta B_{i-2})s\right]e^{-s},\quad \text{for } i=2,3. 
\end{align*}
We define the inner solution as
$U_{inn} := U_0+\varepsilon U_1+\varepsilon^2 U_2+\varepsilon^3 U_3$.
\subsection{Matching}
We first expand the outer solution in the inner variable to $\mathcal{O}(\varepsilon^4)$:
\begin{align*}
    u_{out}(T-\varepsilon s) =& \ u_0(T-\varepsilon s)+\varepsilon u_1(T-\varepsilon s)+\varepsilon^2u_2(T-\varepsilon s)+\varepsilon^3 u_3(T-\varepsilon s) \\
    =&\ \Big(u_0(T)-\varepsilon su_0'(T)+\frac{\varepsilon^2s^2}{2}u_0''(T)-\frac{\varepsilon^3s^3}{6}u_0'''(T)\Big) \\
    &+\varepsilon\Big(u_1(T)-\varepsilon s u_1'(T)+\frac{\varepsilon^2s^2}{2}u_{1}''(T)\Big)
    +\varepsilon^2(u_2(T)-\varepsilon su_{2}'(T))
    +\varepsilon^3u_{3}(T) + \mathcal{O}(\varepsilon^4).
\end{align*}
The inner solution may be written as
\begin{align*}
    u_{inn}(s) =& (A_0+\varepsilon A_1+\varepsilon^2A_2+\varepsilon^3A_3)+s(B_0+\varepsilon B_1+\varepsilon^2 B_2+\varepsilon^3 B_3) \\
    &+\varepsilon^2\left(-\frac{\Delta B_0}{6}s^3+\left(\Delta B_0-\frac{\Delta A_0}{2}\right)s^2+\left[\left(-4\Delta B_0+3\Delta A_0\right)+(\Delta A_0-\Delta B_0)s\right]e^{-s}\right) \\
    &+\varepsilon^3\left(-\frac{\Delta B_1}{6}s^3+\left(\Delta B_1-\frac{\Delta A_1}{2}\right)s^2+\left[\left(-4\Delta B_1+3\Delta A_1\right)+(\Delta A_1-\Delta B_1)s\right]e^{-s}\right).
\end{align*}
Matching these expansions (as $s\rightarrow \infty$), we obtain
\begin{align*}
    &A_0 = u_0(T), A_1 = u_1(T),A_2 = u_2(T),A_3 = u_3(T) \\
    &B_0 = 0, B_1 = -u_0'(T),B_2 = -u_1'(T),B_3= -u_2'(T).
\end{align*}
The PDE problem \eqref{eq: fourth order time PDE given IC} also enforces the polynomial terms in the inner solution to be matched as:
\begin{align*}
    &u_{0}''(T)=\Delta u_0(T) , \qquad 
    u_0'''(T) = \Delta u_0'(T), \qquad
    &u_1''(T)-\Delta u_1(T) =2\Delta u_0'(T)=2u_0'''(T).
\end{align*}
Then, overlap between the inner and outer region is the entire polynomial part of the inner solution. The composite solution to $\mathcal{O}(\varepsilon^4)$ is given by
\begin{align*}
    u_{comp}(x,t) =&\ u_{out}(t)+U_{inn}\left(\frac{T-t}{\varepsilon}\right)-u_{overlap}(t) \\
    =&\ u_0+\varepsilon u_1+\varepsilon^2 u_2+\varepsilon^3 u_3 + \varepsilon^2\left(\Delta u_0(T)+\varepsilon \Delta u_1(T)\right)\left(3+\frac{T-t}{\varepsilon}\right)e^{\frac{t-T}{\varepsilon}} \\
    &+\varepsilon^3\Delta u_0'(T)\left(4+\frac{T-t}{\varepsilon}\right)e^{\frac{t-T}{\varepsilon}}.
\end{align*}
Then, $u^{\varepsilon} = u_{comp}+\mathcal{O}(\varepsilon^4)$. One may directly verify that $u_{comp}$ satisfies $u_{comp,tt}(T) = \mathcal{O}(\varepsilon^2)$ and $u_{comp,ttt}(T) = \mathcal{O}(\varepsilon)$. 

The error rates given in Theorem \ref{thm: convergence} depend on the norms of the higher derivatives of $u^{\varepsilon}$, so it is of interest whether these norms are bounded above independent of $\varepsilon$. The correction terms in $u_{comp}$ coming from the boundary layer will cause these higher derivatives to depend on negative powers of $\varepsilon$. In particular one can directly verify that $\|u_{comp}\|_{H^2(L^2)}$ is bounded above independently of $\varepsilon$. However, each derivative of $u_{comp}$ will pick up an additional factor of $\varepsilon^{-1}$, coming from the exponential term in the boundary layer. For instance, $\|u_{comp}\|_{H^3(L^2)}\sim \varepsilon^{-\frac{1}{2}}$. Crucially, however, for the $\varepsilon$-weighted norm, we compute $\|u_{comp}\|_{H^3_{\varepsilon}(L^2)}\sim e^{-\frac{T}{2\varepsilon}}\varepsilon^{-\frac{1}{2}}$, which may be bounded above independently of $\varepsilon$. This also holds for higher derivatives, that is,
\begin{align*}
    \|u_{comp}\|_{H^k_{\varepsilon}(L^2)}\sim e^{-\frac{T}{2\varepsilon}}\varepsilon^{\frac{5}{2}-k} \sim \left(\frac{2k-5}{Te}\right)^{k-\frac{5}{2}},
\end{align*}
which is independent of $\varepsilon$.
We expect that $u^{\varepsilon}$ also follows this behaviour. This motivates the use the weighted norms in Theorem \ref{thm: convergence} to achieve $\varepsilon$-independent bounds. The higher spatial derivatives are bounded above independent of $\varepsilon$, as the boundary layer only appears in the temporal variable.
\end{document}